\documentclass[a4paper,10pt]{article}

\usepackage[utf8x]{inputenc}
\usepackage{amsfonts}
\usepackage{amsmath}
\usepackage{amssymb}
\usepackage{amsthm}
\usepackage{mathrsfs}
\usepackage{latexsym}
\usepackage{graphicx}
\usepackage{bm}
\usepackage{inputenc}
\usepackage{enumerate}
\usepackage{enumitem}
\usepackage{hyperref}

\swapnumbers

\newtheorem{Lemma}{Lemma}[section]
\newtheorem{Theorem}[Lemma]{Theorem}

\theoremstyle{definition}

\newtheorem{Definition}[Lemma]{Definition}

\theoremstyle{remark}
\newtheorem*{Proof}{Proof}
\newtheorem{Remark}[Lemma]{Remark}

\newtheoremstyle{definition*}
  {3pt}
  {3pt}
  {\rmfamily}
  {}
  { \itshape}
  {.}
  {.5em}
  {\thmnote{#3}}
\theoremstyle{definition*}
\newtheorem*{Definition*}{}

\newtheoremstyle{miniremark*}
  {3pt}
  {3pt}
  {\rmfamily}
  {}
  {\bfseries}
  {.}
  {.5em}
  {\thmnote{#3}}
\theoremstyle{miniremark*}
\newtheorem*{Miniremark*}{}

\DeclareMathOperator{\without}{\sim}
\DeclareMathOperator{\restrict}{\llcorner}

\DeclareMathOperator{\Clos}{Clos}  
\DeclareMathOperator{\Tan}{Tan}     
     
\DeclareMathOperator{\im}{im}       
\DeclareMathOperator{\Lip}{Lip}     
  
\DeclareMathOperator{\dmn}{dmn}     
\DeclareMathOperator{\Nor}{Nor}
    \DeclareMathOperator{\Der}{D}       
\DeclareMathOperator{\pt}{pt}       
\DeclareMathOperator{\ap}{ap}  
\DeclareMathOperator{\loc}{loc}  
\DeclareMathOperator{\gr}{gr} 
\DeclareMathOperator{\reach}{reach}

\newcommand{\Real}[1]{ \mathbf{R}^{#1}}
\newcommand{\Haus}[1]{ \mathscr{H}^{#1} }
\newcommand{\Leb}[1]{ \mathscr{L}^{#1} }
\newcommand{\rect}[1]{(\mathscr{H}^{#1},#1)}
\newcommand{\supHdensity}[3]{\bm{\Theta}^{*#1}(\mathscr{H}^{#1}\restrict \,#2,#3 )}
\newcommand{\Hdensity}[3]{\bm{\Theta}^{#1}(\mathscr{H}^{#1}\restrict \,#2,#3 )}
\newcommand{\supLdensity}[3]{\bm{\Theta}^{*#1}(\mathscr{L}^{#1}\restrict \,#2,#3 )}
\newcommand{\Ldensity}[3]{\bm{\Theta}^{#1}(\mathscr{L}^{#1}\restrict \,#2,#3 )}

\newcommand{\infHdensity}[3]{\bm{\Theta}_{*}^{#1}(\mathscr{H}^{#1}\restrict \,#2,#3 )}
\newcommand{\Lp}[1]{\mathbf{L}_{#1}}

\newcommand{\rsubseteq}{\mathbin{\rotatebox[origin=c]{90}{$\subseteq$}}}

\title{Rectifiability and approximate differentiability of higher order for sets}
\author{Mario Santilli}

\begin{document}
	\maketitle

\begin{abstract}
The main goal of this paper is to develop a concept of approximate differentiability of higher order for subsets of the Euclidean space that allows to characterize higher order rectifiable sets, extending somehow well known facts for functions. We emphasize that for every subset $ A $ of the Euclidean space and for every integer $ k \geq 2 $ we introduce the approximate differential of order $ k $ of $ A $ and we prove it is a Borel map whose domain is a (possibly empty) Borel set. This concept could be helpful to deal with higher order rectifiable sets in applications.
\end{abstract}

\paragraph{\small MSC-classes 2010.}{\small 28A75 (Primary); 49Q15 (Secondary).}

\section{Introduction}	

\textit{For notation and terminology, see the dedicated section in the Introduction.}

\begin{Miniremark*}[Motivation]
	This paper deals with the class of subsets of the $ n $ dimensional Euclidean space that can be covered, up to a set of $ \Haus{m} $ measure zero, by countably many \mbox{$ m $ dimensional} submanifolds of class $ k $, for some integer $ k \geq 1 $. These sets are called \textit{countably $\rect{m}$ rectifiable of class $ k $}; see the definition of \textit{Higher order rectifiability} later in this section.  Evidently, this class of sets finds applications in the study of differentiability properties of ``singular submanifolds'' in geometry and analysis. In what follows, we mention some of these applications.
 
In Calculus of Variations, solutions of problems involving elliptic functionals can be efficiently modeled using the class of \textit{varifolds}, originally introduced by Almgren in the 60's (the classical reference is \cite{MR0307015}). Even just considering the class of \textit{integral varifolds} $ V $ in $ \Real{n} $ whose first variation with respect to area is represented by integrating a function\footnote{See \cite[2.4.12]{MR0257325} for the definition of Lebesgue spaces.} in $ \Lp{\infty}^{\loc}(\|V\|, \Real{n} ) $ (usually called mean curvature), it is well known that there are examples whose support does not locally correspond to a graph of a function of class $ 1 $, on a set of positive $ \Haus{m} $ measure (see \cite[8.1(2)]{MR0307015} for a $ 1 $ dimensional example). However, combining the recent result in \cite{zbMATH06157228} with \cite[8.3]{MR0307015}, we can deduce that the support of an \textit{integral varifold} $ V $ with mean curvature in $ \Lp{m}^{\loc}(\|V\|,\Real{n}) $ is countably $ \rect{m} $ rectifiable of class $ 2 $. Yet, it is unknown if ``class $ 2 $'' in the previous assertion can be replaced by ``class $ k $, for every integer $ k \geq 1 $'', in case the first variation of $ V $ with respect to area equals zero. 

The level sets of maps between Euclidean spaces with distributional derivatives up to order $ k $ representable by integration can be covered, up to a set of $ \Haus{m} $ measure zero, by countably many submanifolds of class $k$ of the appropriate dimension, see \cite[1.6]{MR2183045}.
	
	In Convex Geometry, \cite[Theorem 3]{MR1384392} asserts that if for a convex subset of the $n$ dimensional Euclidean space we define the $ m $-th stratum as the set of points where the normal cone has dimension at least $ n-m$, then the $ m$-th stratum is countably $\rect{m}$ rectifiable of class $ 2 $. This result has been recently extended to every closed subsets of Euclidean space in \cite{2017arXiv170309561M}.
	
Finally we mention the work in \cite{MR2904134} in the context of Legendrian currents and the work in \cite{2015arXiv150600507K} in relation to the notion of discrete curvatures for sets.
\end{Miniremark*}

\begin{Miniremark*}[Results of the present paper]	
	The main contribution of this paper is to introduce a notion of approximate differentiability of higher order for subsets of the Euclidean space and to use it in order to characterize higher order rectifiable sets. For functions whose domain is a subset of the Euclidean space this is a well known fact, established in \cite{MR0043878}, \cite[\S \ 3.1]{MR0257325} and \cite{MR897693}. More specifically these results can be combined to get the following result, see \ref{rectifiability and approximate differentiability for functions I} and \ref{rectifiability and approximate differentiability for functions}.
	\begin{Theorem}[Federer, Isakov, Whitney]
		If $ 1 \leq m < n $ and $ k \geq 1 $ are integers, $ 0 \leq \alpha \leq 1 $, $ A \subseteq \Real{m} $ is $ \Leb{m} $ measurable and $ f : A \rightarrow \Real{n-m} $ is $ \Leb{m} \restrict A $ measurable, then $ f $ is approximately differentiable\footnote{For the definition of approximate differentiability for functions, see \ref{approximate differentiability for functions}.} of order $ (k,\alpha) $ at $ \Leb{m} $ a.e.\ $ a \in A $ if and only if there exist countably many functions $ g_{j}: \Real{m} \rightarrow \Real{n-m} $ of class $ (k,\alpha) $ such that
		\begin{equation*}
		\textstyle \Leb{m}\left( A \sim \bigcup_{j=1}^{\infty}\{ x : g_{j}(x) = f(x)  \} \right) = 0. 
		\end{equation*} 
	\end{Theorem}
	
	In this paper we establish this result for subsets of Euclidean space. In fact, employing the notion of approximate differentiability of higher order for sets introduced in \ref{ap diff for sets} we can prove, in \ref{approximate differentiability of rectifiable sets} and \ref{the set of approximate differentiable points is rectifiable}, the following result.
	
	\begin{Theorem}\label{generalization of Federer-Simon theorem}
		If $ 1 \leq m \leq n $ and $ k \geq 1 $ are integers, $ 0 \leq \alpha \leq 1 $, $ A \subseteq \Real{n} $ is $ \Haus{m} $ measurable and $ \Haus{m}(A) < \infty $, then $ A $ is approximately differentiable of order $ (k,\alpha) $ at $ \Haus{m} $ a.e.\ $ a \in A $ if and only if $ A $ is $ \rect{m} $ rectifiable of class $ (k,\alpha) $.
	\end{Theorem}
	
	It is worth to compare this result with other results in the literature. Firstly, this result can be seen as a generalization to the case of higher order rectifiability of the well known fact in Geometric Measure Theory that $ \rect{m} $ rectifiable sets\footnote{By \cite[3.2.29]{MR0257325} the notion of rectifiability of class $ 1 $ coincides with the classical notion of rectifiability phrased in terms of images of Lipschitzian maps, see \cite[3.2.14]{MR0257325}.} \mbox{of class $ 1 $} can be characterized among all the $ \Haus{m} $ measurable subsets of $ \Real{n} $ with finite $\Haus{m}$ measure through the existence of an $ m $ dimensional ``measure theoretic tangent space'' at $ \Haus{m} $ a.e.\ points of the set. There are essentially two natural ways to define this notion of measure theoretic tangency. One uses a blow up procedure and the other one uses densities of Hausdorff measures.
	
	\begin{Definition}[Simon\footnote{See \cite[11.2, 11.4]{MR756417} and \cite[2.2]{MR1686704}.}]\label{Simon appr tangent plane}
		Suppose $ 1 \leq m \leq n $ are integers, $ A \subseteq \Real{n} $ and \mbox{$ a \in \Real{n} $.} An $ m $ dimensional plane $ T \in \mathbf{G}(n,m) $ is the $ m $ dimensional approximate tangent plane of $ A $ at $ a $ if and only if there exists $ 0 < \theta < \infty $ such that
		\begin{equation*}
		\lim_{r \to 0+}r^{-m} \int_{A}f((x-a)/r)d\Haus{m}x = \theta \int_{T}f d\Haus{m} \quad \textrm{whenever $ f \in \mathscr{K}(\Real{n}) $}.
		\end{equation*}
	\end{Definition}
	
	\begin{Definition}[Federer\footnote{See \cite[3.2.16]{MR0257325}.}]\label{Federer tangent cone}
		Suppose $ 1 \leq m \leq n $ are integers, $ A \subseteq \Real{n} $ and $ a \in \Real{n} $. A vector $ v \in \Real{n} $ is an $ (\Haus{m}\restrict A,m) $ approximate tangent vector at $ a $ if and only if $ \supHdensity{m}{A \cap \mathbf{E}(a,v,\epsilon)}{a} > 0 $ for every $ \epsilon > 0 $. An $ m $ dimensional plane $ T \in \mathbf{G}(n,m) $ is the $ m $ dimensional approximate tangent plane of $ A $ at $ a $ if and only if $ T $ equals the set of all $ (\Haus{m}\restrict A,m) $ approximate tangent vectors at $ a $.
	\end{Definition}
	
	In \ref{Simon appr tangent plane} and in \ref{Federer tangent cone} it is not difficult to see that $ m $ and $ T $ are uniquely determined by $ A $ and $ a $. Either employing the notion in \ref{Simon appr tangent plane} or the one in \ref{Federer tangent cone}, the following well known result holds.
	
\begin{Theorem}[Federer\footnote{See \cite[3.2.19, 3.3.17]{MR0257325}.}, Simon\footnote{See \cite[11.6, 11.8]{MR756417}.}]\label{Federer-Simon theorem}
		Suppose $ 1 \leq m \leq n $ are integers and $ A \subseteq \Real{n} $ is $ \Haus{m} $ measurable with $ \Haus{m}(A) < \infty $. Then $ A $ is $ \rect{m} $ rectifiable of \mbox{class $ 1 $} if and only if $ A $ admits the $ m $ dimensional approximate tangent plane at $ \Haus{m} $ a.e. $ a \in A $.
	\end{Theorem}
	
	Suppose now $ A \subseteq \Real{n} $ and $ a \in \Real{n} $.  The definition of approximate differentiability of order $ 1 $ introduced in \ref{ap diff for sets} (see also \ref{definition of approximate tangent space}) is equivalent to require the existence of a measure theoretic tangent space at $ a $, denoted by $ \ap \Tan(A,a) $. If $ A $ is $ \Haus{m} $ measurable and $ \Haus{m}(A) < \infty $, the sets of points of $ A $ where the aforementioned approximate tangent spaces are $ m $ dimensional subspaces, coincide up to a set of $ \Haus{m} $ measure zero (notice \ref{comparison with Simon approximate tangent plane}, \ref{basic characterization}, \ref{example on tangent cones} and \ref{approximate differentiable set with infinite density}).
	
	The problem of generalizing \ref{Federer-Simon theorem} to the case of higher order rectifiability was addressed in \cite{MR1285779}. In that paper a notion of differentiability of \mbox{order $ 2 $} and order $ (1,\alpha) $ for every $ 0 < \alpha \leq 1 $, is introduced by means of a blow up procedure similar to \ref{Simon appr tangent plane}. However, as it is pointed out in \cite[pp.\ 7--8]{MR1285779}, examples show that $ \rect{m} $ rectifiable sets of class $ 2 $ may fail to be differentiable of \mbox{order $ 2 $} in the sense of \cite{MR1285779} at $ \Haus{m} $ a.e.\ points. Therefore, in order to generalize \ref{Federer-Simon theorem}, additional technical hypotheses on the structure of the sets are needed in the main theorems \cite[3.5, 3.12]{MR1285779}. These facts suggest the possibility to consider a different notion of (approximate) differentiability of order greater than $ 1 $ and in the present paper we accomplish such goal. Our notion is based on the approach of \ref{Federer tangent cone} (see also \ref{Mattila tangent plane}), rather than \ref{Simon appr tangent plane}.
	
	For every integer $ k \geq 2 $, the notion of approximate differentiability of \mbox{order $ k $} for a subset $ A \subseteq \Real{n} $ naturally induces a notion of approximate differential of order $ k $, $ \ap\Der^{k} A $, of $ A $; see \ref{definition of approximate differentials}. For every $ A \subseteq \Real{n} $, this is always a Borel map with values in $ \textstyle\bigodot^{k}(\Real{n},\Real{n}) $ whose domain is a (possibly empty) Borel subset of $ \Real{n} $, see \ref{approx differentials are Borel maps}. Moreover the approximate differential of order $ 2 $ naturally induces a notion of ``approximate second fundamental form''. In fact, for every $ a \in \dmn \ap \Der^{2}A $ this can be defined as the symmetric bilinear form
	\begin{equation*}
	\ap \Der^{2}A(a) | \ap \Tan(A,a) \times \ap \Tan(A,a).
	\end{equation*}
	In \ref{rectifiable sets and normal vector fields} and \ref{touching by balls and estimate ssf} two classical properties of the second fundamental form of submanifolds of class $ 2 $ are extended to our setting.
	
	Finally we mention that a notion of pointwise differentiability for subsets of the Euclidean space has been recently developed in \cite{snulmenn:sets.v1} to study higher order differentiability properties of stationary varifolds. In the present paper we establish the connection between the notion of approximate differentiability and pointwise differentiability in \ref{ap diff and pt diff}. 
\end{Miniremark*}

\begin{Miniremark*}[Notation and basic definitions] 
The notation and the terminology used withouth comments agree with \cite[pp.\ 669-676]{MR0257325}. However, for the reader's convenience, sometimes we use footnotes to point out the references \mbox{in \cite{MR0257325}.} Moreover we add the following classical definitions.

\begin{Definition*}[Distance function from a set]
 Let $ A \subseteq \Real{n} $. We define $ \bm{\delta}_{A} $ to be the function on $ \Real{n} $ such that $ \bm{\delta}_{A}(x) = \inf \{ |x-a| : a \in A \} $ for $ x \in \Real{n} $.
\end{Definition*}

\begin{Definition*}[Orthogonal projections]
If $ 1 \leq m \leq n $ are integers we define $ \mathbf{G}(n,m) $ to be the set of all $ m $ dimensional subspaces of $ \Real{n} $. If $ T \in \mathbf{G}(n,m) $ we define \mbox{$ T_{\natural} : \Real{n} \rightarrow \Real{n} $} to be the linear map such that
\begin{equation*}
 T_{\natural}^{*} = T_{\natural}, \;\; T_{\natural}\circ T_{\natural} = T_{\natural}, \;\; \im T_{\natural} = T,
\end{equation*}
and we define $ T^{\perp} = \ker T_{\natural} $.
\end{Definition*}

\begin{Definition*}[Pointwise differentiability for functions]
Suppose $ X $ and $ Y $ are normed vector spaces, $ k \geq 0 $ is an integer, $ 0 \leq \alpha \leq 1 $, $ g $ maps a subset of $ X $ into $ Y $ and $ a \in X $. We say that $ g $ is pointwise differentiable of order $ (k,\alpha) $ at $ a $ if and only if there exists an open set $ U \subseteq X $ and a polynomial function $ P : X \rightarrow Y $ of degree at most $ k $ such that $ a \in U \subseteq \dmn g $, $ g(a) = P(a) $,
\begin{equation*}
 \lim_{x\to a} \frac{|g(x) - P(x)|}{|x-a|^{k}} = 0 \;\;\textrm{if $ \alpha = 0 $,}\quad \limsup_{x \to a}\frac{|g(x)-P(x)|}{|x-a|^{k + \alpha}}< \infty \;\; \textrm{if $ \alpha >0 $}.
\end{equation*}
In this case $ P $ is unique and the pointwise differentials of order $ i $ of $ f $ at $ a $ are defined by $ \pt \Der^{i}g(a) = \Der^{i}P(a) $ for $ i = 0, \ldots , k $.
\end{Definition*}

\begin{Definition*}[Functions and submanifolds of class $ (k,\alpha) $]
Suppose $ X $ and $ Y $ are normed vector spaces, $ k \geq 0 $ is an integer, $ 0 \leq \alpha \leq 1 $, $ g $ maps some open subset of $ X $ into $ Y $ and $ a \in X $. We say that $ g $ is of class $ (k,\alpha) $ if and only if $ g $ is of \mbox{class $ k $} and each point of $ \dmn g $ has an open neighbourhood $ U $ such that $ (\Der^{k}f)|U $ satisfies a H\"older condition with exponent $ \alpha $.

Suppose $ k \geq 0 $ is an integer and $ 0 \leq \alpha \leq 1 $. The notion of diffeomorphism of class $ (k,\alpha) $ is made by replacing ``class $ k $'' with ``class $ (k,\alpha) $'' in \cite[3.1.18]{MR0257325}. Analogously the notion of $ \mu $ dimensional submanifold of class $ (k,\alpha) $ of $ \Real{n} $ is made by replacing ``class $ k $'' with ``class $ (k,\alpha) $'' in \cite[3.1.19]{MR0257325}.
\end{Definition*}

\begin{Definition*}[Second fundamental form]
If $ 1 \leq m \leq n $ are integers, $ M $ is an $ m $ dimensional submanifold of class $ 2 $ of $ \Real{n} $ and $ a \in M $ then we call second fundamental form of $ M $ at $ a $ the unique symmetric $ 2 $ linear function
\begin{equation*}
 \mathbf{b}_{M}(a) : \Tan(M,a) \times \Tan(M,a) \rightarrow \Nor(M,a)
\end{equation*}
such that $ \mathbf{b}_{M}(a)(u,v) \bullet \nu(a) = - \Der \nu(a)(u) \bullet v $ for each $ u, v \in \Tan(M,a) $, whenever $ \nu : M \rightarrow \Real{n} $ is of class $ 1 $ relative to $ M $ with $ \nu(x) \in \Nor(M,x) $ for every $ x \in M $.
\end{Definition*}

\begin{Definition*}[Cones]
		A subset $ C \subseteq \Real{n} $ is called cone if and only if $ \lambda x \in C $ whenever $ x \in C $ and $ \lambda > 0 $. 
	\end{Definition*}

\begin{Definition*}[Higher order rectifiability\footnote{When $ \phi = \Haus{m} $ this notion has been introduced in \cite[3.1]{MR1285779}.}]
		Suppose $ 1 \leq m \leq n $ are integers and $ \phi $ is a measure over $ \Real{n} $.  A subset $ A \subseteq \Real{n} $ is called \textit{countably $ (\phi,m) $ rectifiable of class $ (k,\alpha) $} if and only if there exist countably many $ m $ dimensional submanifolds $ M_{j} $ of class $ (k,\alpha) $ such that 
		\begin{equation*}
			\textstyle \phi\big( A \sim \bigcup_{j=1}^{\infty} M_{j} \big) =0.
		\end{equation*}
		A subset $ A \subseteq \Real{n} $ is called \textit{$ (\phi,m) $ rectifiable of class $ (k,\alpha) $} if it is countably $ (\phi,m) $ rectifiable of class $ (k,\alpha) $ and $ \phi(A) < \infty $.
	\end{Definition*}

\begin{Definition*}[Some further notation]
	Suppose $ 1 \leq m \leq n $ and $ k \geq 1 $ are integers, $ 0 \leq \alpha \leq 1 $, $ a \in \Real{n} $, $ T \in \mathbf{G}(n,m) $, $ 0 \leq \kappa < \infty $ and suppose $ f : T \rightarrow T^{\perp} $ is a function such that $ f(T_{\natural}(a))=T^{\perp}_{\natural}(a) $.

Then we define\footnote{Compare this definition with similar ones introduced in \cite[3.3.1]{MR0257325} and \cite[15.12]{MR1333890}.}
	\begin{equation*}
		\mathbf{X}_{k,\alpha}(a,T,f,\kappa) = \Real{n} \cap \{ z : | f(T_{\natural}(z)) - T^{\perp}_{\natural}(z) | \leq \kappa |T_{\natural}(z-a)|^{k+\alpha}    \};
	\end{equation*}
alternatively $ \mathbf{X}_{k}(a,T,f,\kappa) = \mathbf{X}_{k,0}(a,T,f,\kappa) $. If $ f(\chi) = T^{\perp}_{\natural}(a) $ for every $ \chi \in T $ then we abbreviate $ \mathbf{X}(a,T,\kappa) = \mathbf{X}_{1}(a,T,f,\kappa) $.
	
If $ 0 < s < \infty $ and $ 0 < t < \infty $ we define
\begin{equation*}
 \mathbf{C}(T,a,s,t) = \Real{n} \cap \{ x : |T_{\natural}(x-a)| < s, \; |T^{\perp}_{\natural}(x-a)| < t \}.
\end{equation*}

Finally let $ \gr(f) = \{ \chi + f(\chi) : \chi \in T \} $.
\end{Definition*}

\end{Miniremark*}

\begin{Miniremark*}[Organization of the paper]
In section \ref{section :Approximate differentiability for functions} we recall the theory of approximate differentiability for functions because both we use it in the following sections and it provides a scheme for the theory of approximate differentiability for sets we develop later. In section \ref{section: Approximate differentiability for sets} the key concepts of lower and upper approximate tangent cones (see \ref{approx tangent cones}), approximate differentiability (see \ref{ap diff for sets}) and approximate differentials (see \ref{definition of approximate differentials}) are introduced together with proofs of a basic characterization in \ref{basic characterization}, one part of \ref{generalization of Federer-Simon theorem} in \ref{approximate differentiability of rectifiable sets} and some illustrative examples in \ref{example on tangent cones} and \ref{approximate differentiable set with infinite density}. At the end of section \ref{section: Approximate differentiability for sets} we generalize to rectifiable sets of class $ 2 $ the classical equation relating the differential of a normal vector field and the second fundamental form of a submanifold of \mbox{class $ 2 $.} In section \ref{section Relation with pointwise differentiability}, after giving in \ref{basic remark on pt tangent cone} an equivalent formulation of the notion of pointwise differentiability of \mbox{order $ 1 $} for sets, we prove a result enlightening the relation between approximate differentiability and pointwise differentiability of higher order for sets in \ref{ap diff and pt diff} and a basic estimate for the approximate second fundamental form in \ref{touching by balls and estimate ssf}. Finally in section \ref{section: Rectifiability and Borel measurability} we prove that the approximate differentials are Borel maps and the remaining part of \ref{generalization of Federer-Simon theorem}.

The content of this paper was part of author's PHD thesis, supervised by Ulrich Menne, submitted at the University of Potsdam.
\end{Miniremark*}

\begin{Miniremark*}[Acknowledgements]
The author is grateful to Ulrich Menne, who suggested this problem, carefully read the original manuscript and provided the author with a very detailed list of comments, corrections and improvements. 

This work was developed while the author was financially supported by the ``IMPRS for Geometric Analysis, Gravitation and String Theory'' and the ``IMPRS for Mathematical and Physical Aspects for Gravitation, Cosmology and Quantum Field Theory''.
\end{Miniremark*}

\section{Approximate differentiability for functions} \label{section :Approximate differentiability for functions}

\begin{Definition}
	Let $ f $ be a function mapping a subset of $ \Real{m} $ into some set $ Y $ and let $ a \in \Real{m} $. If $ Y $ is a normed vector space, a point $ y \in Y $ is the approximate limit of $ f $ at $ a $ if and only if
	\begin{equation*}
	\Ldensity{m}{\Real{m} \sim \{ x : | f(x) - y |  \leq \epsilon \}}{a} = 0 \quad \textrm{for every $ \epsilon > 0 $}
	\end{equation*}
	and we denote it by $ \ap\lim_{x \to a} f(x) $. If $ Y = \overline{\Real{}} $, a point $ t \in \overline{\Real{}} $ is the approximate upper limit of $ f $ at $ a $ if and only if
	\begin{equation*}
	t = \inf\{ s: \Ldensity{m}{ \{ x :  f(x) > s \}}{a} = 0 \}
	\end{equation*}
	and we denote it by $ \ap \limsup_{x \to a} f(x) $.
\end{Definition}

\begin{Remark}
	This concept is a special case of \cite[2.9.12]{MR0257325}.
\end{Remark}

\begin{Definition}\label{approximate differentiability for functions}
	Let $ 1 \leq m < n $ and $ k \geq 0 $ be integers, $ 0 \leq \alpha \leq 1 $, $ A \subset \mathbf{R}^{m} $, $ f : A \rightarrow \Real{n-m} $ and $ a \in \mathbf{R}^{m} $.
	
	We say that $ f $ is \textit{approximately differentiable of order $ (k,\alpha) $ at $ a $} (\textit{$ f $ is approximately differentiable of order $ k $ at $ a $} if $ \alpha =0 $) if
	\begin{equation*}
		\Ldensity{m}{\mathbf{R}^{m}\without A}{a} =0
	\end{equation*}
	and there exists a polynomial function $ P :\Real{m} \rightarrow \Real{n-m} $ of degree at most $ k $ such that $ P(a)= f(a) $ if $ a \in A $,
	\begin{equation*}
		\ap\lim_{x \to a} \frac{|f(x)-P(x)|}{|x-a|^{k}}=0 \;\; \textrm{if $ \alpha =0 $}, \quad \ap\limsup_{x \to a}\frac{|f(x)-P(x)|}{|x-a|^{k + \alpha}}< \infty  \;\; \textrm{if $ \alpha >0 $}.
	\end{equation*}
\end{Definition}

\begin{Remark}
	The condition $\Ldensity{m}{\mathbf{R}^{m}\without A}{a} =0 $ in \ref{approximate differentiability for functions} is redundant if \mbox{$ \alpha =0 $.} Moreover, employing a classical result due to De Giorgi, see \cite[Lemma 2.I]{MR0167862}, we deduce that the polynomial function $ P $ in \ref{approximate differentiability for functions} is uniquely determined by $ f $ and $ a $.
\end{Remark}

\begin{Definition}
	Let $ A \subset \Real{m} $ and let $ f : A \rightarrow \Real{n-m} $. For every non negative integer $ k $ the function $ \ap\Der^{k} f $ is defined to be the function whose domain consists of all $ a \in \Real{m} $ such \mbox{that $ f $} is approximately differentiable of order $ k $ at $ a $ and whose value at $ a $ equals $ \Der^{k}P(a) $, where $ P $ satisfies \ref{approximate differentiability for functions}.
\end{Definition}

\begin{Remark}
	If $ a \in A \subset \Real{m} $ and $ f : A \rightarrow \Real{n-m} $ then $ f $ is approximately differentiable of order $ 0 $ at $ a $ if and only if $ f $ is $ (\Leb{m},V) $ approximately continuous\footnote{See \cite[2.9.12]{MR0257325}.} \mbox{at $ a $}. In this case $ \ap\Der^{0}f(a)=f(a) $. Here $ V $ is the standard $ \Leb{m} $ Vitali relation, $ V = \{ (a, \mathbf{B}(a,r)): a \in \Real{m}, \; 0 < r < \infty \} $.
	
	In case $ a \in A $ the notion of approximate differentiability of order $ 1 $ has been introduced in \cite[3.1.2]{MR0257325}.
\end{Remark}

\begin{Lemma}\label{alternative characterization of approximate differentiability for functions}
	Suppose $ 1 \leq m < n $ are integers, $ A \subseteq \Real{m} $, $ a \in \Real{m} $, $ f : A \rightarrow \Real{n-m} $, $ \gamma \geq 1 $, $ 0 < M < \infty $ and $ 0 \leq \lambda < \infty $  such that
	\begin{equation*}
		\limsup_{r \to 0+} \frac{\Leb{m}(\mathbf{B}(a,r) \cap \{ x : |f(x)| > \lambda \,r^{\gamma}\})  }{ \bm{\alpha}(m)\,r^{m}} < M.
	\end{equation*}
	
	Then $ \supLdensity{m}{ \{ x : |f(x)| > 2^{\gamma}\lambda\,|x-a|^{\gamma} \} }{a} < M(1-2^{-m})^{-1} $.
\end{Lemma}

\begin{Proof}
	Let $ \delta > 0 $ such that 
	\begin{equation*}
		\Leb{m}(\mathbf{B}(a,r) \cap \{ x : |f(x)| > \lambda \,r^{\gamma}\}) < M\,\bm{\alpha}(m)\,r^{m} \quad \textrm{for $ 0 < r \leq \delta $}.
	\end{equation*}
	Therefore for $ 0 < r \leq \delta $ we observe
	\begin{flalign*}
		& \mathbf{B}(a,r) \cap \{  x : |f(x)| > 2^{\gamma} \lambda \,|x-a|^{\gamma}\} \\
		& \quad = \{a\} \cup \bigcup_{i=0}^{\infty} (\mathbf{B}(a,r/2^{i}) \sim \mathbf{B}(a,r/2^{i+1})) \cap \{  x : |f(x)| > 2^{\gamma} \lambda \,|x-a|^{\gamma}\} \\
		& \quad \subseteq \{a\} \cup \bigcup_{i=0}^{\infty} \mathbf{B}(a,r/2^{i}) \cap \{ x : |f(x)|> \lambda\,(r/2^{i})^{\gamma}\},
	\end{flalign*}
	\begin{equation*}
		\Leb{m}(\mathbf{B}(a,r) \cap \{  x : |f(x)| > 2^{\gamma} \lambda \,|x-a|^{\gamma}\}) < M\,\bm{\alpha}(m)r^{m}(1-2^{-m})^{-1}
	\end{equation*}
	and the conclusion follows.
\end{Proof}

\begin{Theorem}\label{alternative definition for approximate differentiability for functions}
	Let $ 1 \leq m < n $ and $ k \geq 1 $ be integers, $ 0 \leq \alpha \leq 1 $, $ A \subset \Real{m} $, $ a \in \Real{m} $ and $ f : A \rightarrow \Real{n-m} $. 
	
	Then $ f $ is approximately differentiable of order $ (k,\alpha) $ at $ a $ if and only if there exists a function $ g : \Real{m} \rightarrow \Real{n-m} $ pointwise differentiable of order $ (k,\alpha) $ at $ a $ such that $ f(a) = g(a) $ if $ a \in A $ and
	\begin{equation*}
		\Ldensity{m}{ \Real{m} \without \{  x: g(x)=f(x)\}  }{a} =0.
	\end{equation*}
	In this case $ \pt\Der^{i}g(a) = \ap \Der^{i}f(a) $ for $ i= 0, \ldots , k $.
\end{Theorem}

\begin{Proof}
	Suppose $ f $ is approximately differentiable of order $ (k,\alpha) $ at $ a $ and $ \alpha = 0 $. There exists a polynomial function \mbox{$ P : \Real{m} \rightarrow \Real{n-m} $} of degree at most $ k $ such that, if for every integer $ i \geq 1 $ we define $ 	S_{i}= \{  x : | f(x)-P(x) | < i^{-1}|x-a|^{k} \} $, then there exists $ \delta_{i}>0 $ such that $ \Leb{m}(\mathbf{B}(a,r) \sim S_{i}) < 2^{-i}r^{m} $ for $ 0 < r \leq \delta_{i} $. We can assume $ \delta_{i+1} < \delta_{i} $ for each $ i \geq 1 $ and $ \delta_{i} \to 0 $ as $ i \to \infty $. Let
	\begin{equation*}
		\textstyle T = \bigcup _{i=1}^{\infty}\left[ S_{i} \cap \mathbf{B}(a,\delta_{i}) \sim \mathbf{B}(a,\delta_{i+1})\right].
	\end{equation*}
	If $ r > 0 $  and $ j \geq 1 $ is an integer such that $ \delta_{j+1} < r \leq \delta_{j} $ we compute
	\begin{equation*}
		\Leb{m}(\mathbf{B}(a,r) \sim T) \leq \Leb{m}(\mathbf{B}(a,r) \sim S_{j}) + \sum_{l=j+1}^{\infty}\Leb{m}( \mathbf{B}(a,\delta_{l}) \sim S_{l}) < r^{m} \sum_{l=j}^{\infty} 2^{-l}
	\end{equation*}
	and we conclude $ \Ldensity{m}{\Real{m} \sim T}{a} =0 $. Moreover
	\begin{equation*}
		\lim_{T \ni x \to a} \frac{| f(x) - P(x)|}{|x-a|^{k}}=0.
	\end{equation*}
	If we define $ g : \Real{m} \rightarrow \Real{n-m} $ as $ g(x)=f(x) $ if $ x \in T $ and $ g(x)= P(x) $ if $ x \in \Real{m} \sim T $, then we have $ \Ldensity{m}{\Real{m} \sim \{  x: g(x)=f(x)\} }{a}=0 $,
	\begin{equation*}
		\lim_{ x \to a} \frac{| g(x) - P(x)|}{|x-a|^{k}}=0 \quad \textrm{and} \quad g(a)= P(a),
	\end{equation*}
since $ a \notin T $. If $ \alpha > 0 $, once we have chosen $ 0 \leq \lambda < \infty $ such that 
	\begin{equation*}
		\ap\limsup_{x \to a}\frac{|f(x)-P(x)|}{|x-a|^{k + \alpha}}< \lambda,
	\end{equation*}
	we can use the same argument above replacing the sets $ S_{i} $ with the set
	\begin{equation*}
		S = \{  x : | f(x)-P(x) | < \lambda\,|x-a|^{k+ \alpha} \}.
	\end{equation*}
\end{Proof}

\begin{Remark}
	The proof of \ref{alternative definition for approximate differentiability for functions} has been adapted from \cite[3.2.16]{MR0257325} and \cite[3.1.22]{MR0257325}.
\end{Remark}

\begin{Remark}\label{approximate differentiability of rectifiable functions}
	Let $ 1 \leq m < n $ and $ k \geq 1 $ be integers, $ 0 \leq \alpha \leq 1 $, let $ A \subset \Real{m} $ be $ \Leb{m} $ measurable and let $ f : A \rightarrow \Real{n-m} $ be $ \Leb{m}\restrict A $ measurable. If there exist countably many functions $ g_{j} : \Real{m} \rightarrow \Real{n-m} $ of class $ (k,\alpha) $ such that 
	\begin{equation*}
		\textstyle \Leb{m}\left( A \sim \bigcup_{j=1}^{\infty} A_{j}  \right) =0
	\end{equation*}
	where $ A_{j} = A \cap \{ z : g_{j}(z)= f(z)  \} $ for $ j \geq 1 $, then using \cite[2.10.19(4)]{MR0257325} and \ref{alternative definition for approximate differentiability for functions} we can easily prove that $ f $ is approximately differentiable of order $ (k,\alpha) $ at $ \Leb{m} $ a.e.\ $ a \in A $ and, for each $ j \geq 1 $, 
	\begin{equation*}
		\Der^{i}g_{j}(z)= \ap\Der ^{i}f(z) \quad \textrm{for $ \Leb{m} $ a.e.\ $ z \in A_{j} $ and $ i = 0 , \ldots , k $}.
	\end{equation*}
	
\end{Remark}

\begin{Theorem}\label{rectifiability and approximate differentiability for functions I}
	Let $ 1 \leq m < n $ and $ k \geq 0 $ be integers, $ A \subset \Real{m} $ and let $ f : A \rightarrow \Real{n-m} $ be approximately differentiable of order $ (k,1) $ at $ \Leb{m} $ a.e.\ $ a \in A $.
	Then the following statements hold.
	\begin{enumerate}
		\item \label{rectifiability and approximate differentiability for functions I: rademacher} $ f $ is approximately differentiable of order $ k+1 $ at $ \Leb{m} $ a.e. $ x \in A $.
		\item \label{rectifiability and approximate differentiability for functions I: measurability}   $ A $ is $ \Leb{m} $ measurable and the functions $ \ap \Der^{i}f $ are $ \Leb{m}\restrict A $ measurable for $ i = 0, \ldots , k+1 $.
		\item \label{rectifiability and approximate differentiability for functions I: rectifiability} There exist countably many functions $ g_{j}: \Real{m} \rightarrow \Real{n-m} $ of class $ k+1 $ such that
		\begin{equation*}
			\textstyle \Leb{m}\left( A \sim \bigcup_{j=1}^{\infty}\{ x : g_{j}(x) = f(x)  \} \right) = 0. 
		\end{equation*}
	\end{enumerate}
\end{Theorem}

\begin{Proof}
	First we observe that $ A $ is $ \Leb{m} $ measurable, $ f $ is $ (\Leb{m},V) $ approximately continuous\footnote{As usual, $ V = \{ (a, \mathbf{B}(a,r)): a \in \Real{m}, \; 0 < r < \infty \} $.} at $ \Leb{m} $ a.e.\ $ a \in A $ and $ f $ is $ \Leb{m} \restrict A $ measurable by \cite[2.9.11, 2.9.13]{MR0257325}.
	
	If $ k = 0 $ the conclusions are consequences of \cite[3.1.8, 3.1.4, 3.1.16]{MR0257325} respectively. We use induction over $ k $. Since $ f $ is approximately differentiable of \mbox{order $ (k-1,1) $} at $ \Leb{m} $ a.e.\ point of $ A $ we inductively assume that $ \ap \Der^{i}f $ are $ \Leb{m} \restrict A $ measurable for $ i = 0, \ldots , k $. We use now \cite[Theorem 2]{MR897693} and \cite[3.1.15]{MR0257325} to deduce the existence of countably many functions $ g_{j}: \Real{m} \rightarrow \Real{n} $ of class $ k+1 $ satisfying \eqref{rectifiability and approximate differentiability for functions I: rectifiability}. Now \eqref{rectifiability and approximate differentiability for functions I: rademacher} and \eqref{rectifiability and approximate differentiability for functions I: measurability} follow from \ref{approximate differentiability of rectifiable functions}. 
	
\end{Proof}

\begin{Theorem}\label{rectifiability and approximate differentiability for functions}
	Suppose $ 1 \leq m < n $ and $ k \geq 1 $ are integers, $ 0 \leq \alpha \leq 1 $, \mbox{$ A \subset \Real{m} $} and $ f : A \rightarrow \Real{n-m} $ is approximately differentiable of order $ (k,\alpha) $ at $ \Leb{m} $ a.e.\ $ a \in  A $.
	
	Then the following statements hold.
	\begin{enumerate}
		\item \label{rectifiability and approximate differentiability for functions: measurability}  $ A $ is $ \Leb{m} $ measurable and the functions $ \ap \Der^{i}f $ are $ \Leb{m}\restrict A $ measurable for $ i = 0, \ldots, k $.
		\item \label{rectifiability and approximate differentiability for functions: rectifiability} There exist countably many functions $ g_{j}: \Real{m} \rightarrow \Real{n-m} $ of class $ (k,\alpha) $ such that
		\begin{equation*}
			\textstyle \Leb{m}\left( A \sim \bigcup_{j=1}^{\infty}\{ x : g_{j}(x) = f(x)  \} \right) = 0. 
		\end{equation*}
	\end{enumerate}
\end{Theorem} 

\begin{Proof}
	
	Since $ f $ is approximately differentiable of order $ (k-1,1) $ at every $ x \in A $ then \eqref{rectifiability and approximate differentiability for functions: measurability} follows from \ref{rectifiability and approximate differentiability for functions I}\eqref{rectifiability and approximate differentiability for functions I: measurability}. Now we can apply \cite[Theorem 1]{MR897693} if $ \alpha = 0 $ or \cite[Theorem 2]{MR897693} if $ \alpha > 0 $ to get \eqref{rectifiability and approximate differentiability for functions: rectifiability}.
\end{Proof}

\section{Approximate differentiability for sets}\label{section: Approximate differentiability for sets}

\begin{Definition}\label{approx tangent cones}
	Suppose $ X $ is a normed vector space, $ \phi $ is a measure over $ X $, $ m $ is a positive integer and $ a \in X $. 
	
	We define \textit{the $ m $ dimensional approximate upper tangent cone of $ \phi $ at $ a $} by\footnote{As in \cite[3.2.16]{MR0257325}, $ \mathbf{E}(a,v,\epsilon) = X \cap \{x : |r(x-a) - v| < \epsilon \; \textrm{for some $ 0 < r < \infty $} \} $.}
	\begin{equation*}
		\Tan^{*m}(\phi,a) = X \cap \{ v : \bm{\Theta}^{*m}(\phi \restrict \mathbf{E}(a,v,\epsilon),a) > 0 \; \textrm{for every $ \epsilon > 0 $} \}
	\end{equation*}
	and \textit{the $ m $ dimensional approximate lower tangent cone of $ \phi $ at $ a $} as the set $ \Tan^{m}_{*}(\phi,a)$ of $ v \in X $ such that for every $ \epsilon > 0 $ there exists $ \eta > 0 $ such that 
	\begin{equation*}
		\phi( \mathbf{U}(a+rv,\epsilon r)) \geq \eta r^{m} \quad \textrm{whenever $ 0 < r \leq \eta $.}
	\end{equation*}
	
	In case $ \Tan^{*m}(\phi,a) = \Tan_{*}^{m}(\phi,a) $, this set is denoted by $ \Tan^{m}(\phi,a) $ and we call it \textit{the $ m $ dimensional approximate tangent cone of $ \phi $ at $ a $}.
\end{Definition}

\begin{Remark}\label{Remark on lower and upper densities and cones}
	Evidently $ \Tan^{m}_{*}(\phi,a) \subseteq \Tan^{*m}(\phi,a) $. Moreover one may easily verify that $ \Tan_{*}^{m}(\phi,a) $ and $ \Tan^{*m}(\phi,a) $ are closed cones. Finally  
	\begin{equation*}
		\bm{\Theta}^{*m}(\phi,a) > 0 \quad [\bm{\Theta}^{m}_{*}(\phi,a) > 0 ] \iff 0 \in \Tan^{*m}(\phi,a) \quad [0 \in \Tan^{m}_{*}(\phi,a)].
	\end{equation*}
	
\end{Remark}

\begin{Remark}\label{basic properties of ap tangent cones}
	Observe that, in this case, our notation does not agree with \cite[3.2.16]{MR0257325}. In fact, $ \Tan^{*m}(\phi,a) $ is denoted by $ \Tan^{m}(\phi,a) $ in \cite[3.2.16]{MR0257325}.
	
	It is often useful to recall that if $ C $ is a compact subset of $ X \sim \Tan^{*m}(\phi,a) $ and $ T = \{ a + rv : r \geq 0, \, v \in C \} $ then $ \bm{\Theta}^{m}(\phi \restrict T,a) = 0 $. This is proved in \cite[3.2.16]{MR0257325}.
\end{Remark}

\begin{Remark}\label{example on tangent cones}
	It is natural to consider the following cone
	\begin{equation*}
		T = X \cap \{ v : \bm{\Theta}^{m}_{*}(\phi \restrict \mathbf{E}(a,v,\epsilon),a) > 0 \; \textrm{for every $ \epsilon > 0 $} \}.
	\end{equation*}
	Evidently $ \Tan^{m}_{*}(\phi,a) \subseteq T $, but simple examples show that the opposite inclusion does not hold. In fact, if we consider $ X = \Real{} $, $ \phi = \Leb{1} \restrict A $, $ m = 1 $ and $ a = 0 $, where 
	\begin{equation*}
		A =  \bigcup_{i=0}^{\infty} \Real{} \cap \{ t : 2^{-2i-1} < |t| < 2^{-2i} \},
	\end{equation*}
	then $ \Tan^{1}_{*}(\phi,0) = \{0\} $ and $ T = \Real{} $.
\end{Remark}

\begin{Remark}\label{Density and inclusion of approximate tangent cones}
	Suppose $ 1 \leq m \leq n $ are integers, $ A \subseteq \Real{n} $, $ B \subseteq \Real{n} $ and $ a \in \Real{n} $. 
	
	If $ \Hdensity{m}{A \sim B}{a} = 0 $ then it is not difficult to see that
	\begin{equation*}
		\Tan^{m}_{*}(\Haus{m}\restrict A,a) \subseteq \Tan^{m}_{*}(\Haus{m}\restrict B , a),
	\end{equation*}
	\begin{equation*}
		\Tan^{*m}(\Haus{m}\restrict A,a) \subseteq \Tan^{*m}(\Haus{m}\restrict B,a).
	\end{equation*}
\end{Remark}

\begin{Lemma}\label{remark on approximate tangent cone I}
	Suppose $ 1 \leq m \leq n $ are integers, $ A \subseteq \Real{n} $, $ a \in \Real{n} $ and \mbox{$ T \in \mathbf{G}(n,m) $.}
	
	Then the following three conditions are equivalent:
	\begin{enumerate}
		\item \label{remark on approximate tangent cone I:inclusion} $ \Tan^{*m}(\Haus{m}\restrict A,a) \subseteq T $,
		\item \label{remark on approximate tangent cone I:cone condition} $ \Hdensity{m}{A \sim \mathbf{X}(a,T,\epsilon)}{a} =0 $ whenever $ \epsilon > 0 $,
		\item \label{remark on approximate tangent cone I:vertical dist condition} whenever $ \epsilon > 0 $
		\begin{equation*}
			\lim_{r \to 0}\frac{\Haus{m}( A \cap \mathbf{B}(a,r) \cap \{ z : |T^{\perp}_{\natural}(z-a)| > \epsilon\,r\})  }{\bm{\alpha}(m)\,r^{m}}=0.
		\end{equation*}
	\end{enumerate}
\end{Lemma} 
\begin{Proof}
	The fact that \eqref{remark on approximate tangent cone I:inclusion} implies \eqref{remark on approximate tangent cone I:cone condition} is a consequence of \ref{basic properties of ap tangent cones} and the fact that \eqref{remark on approximate tangent cone I:vertical dist condition} follows from \eqref{remark on approximate tangent cone I:cone condition} is evident. If the condition in \eqref{remark on approximate tangent cone I:vertical dist condition} holds for some $ \epsilon > 0 $ then we can argue as in \ref{alternative characterization of approximate differentiability for functions} to show that 
	\begin{equation*}
		\lim_{r \to 0}\frac{\Haus{m}( A \cap \mathbf{B}(a,r) \cap \{ z : |T^{\perp}_{\natural}(z-a)| > 2 \epsilon\,|z-a| \})  }{\bm{\alpha}(m)\,r^{m}}=0.
	\end{equation*}
	Therefore \eqref{remark on approximate tangent cone I:vertical dist condition} implies \eqref{remark on approximate tangent cone I:inclusion}.
\end{Proof}

\begin{Lemma}\label{alternative characterization of approximate differentiability for sets}
	Let $ 1 \leq m \leq n $ and $ k \geq 1 $ be integers, $ 0 \leq \alpha \leq 1 $, $ 0 \leq \lambda < \infty $, $ 0 < M < \infty $, $ A \subseteq \Real{n} $, $ a \in \Real{n} $, $ T \in \mathbf{G}(n,m) $ and let $ f : T \rightarrow T^{\perp} $ be a function of class $ 1 $ such that $ f(T_{\natural}(a)) = T^{\perp}_{\natural}(a) $ and $ \Der f(T_{\natural}(a)) = 0 $. Suppose 
	\begin{equation*}
		\lim_{r \to 0}\frac{ \Haus{m}( A \cap \mathbf{B}(a,r) \cap \{ z: |T^{\perp}_{\natural}(z) - f(T_{\natural}(z))| > \epsilon\,r \})  }{\bm{\alpha}(m)\,r^{m}} = 0 \; \textrm{for every $ \epsilon > 0 $},
	\end{equation*}
	\begin{equation*}
		\limsup_{r \to 0}\frac{ \Haus{m}( A \cap \mathbf{B}(a,r) \cap \{ z: |T^{\perp}_{\natural}(z) - f(T_{\natural}(z))| > \lambda\,r^{k + \alpha}\})  }{\bm{\alpha}(m)\,r^{m}} < M.
	\end{equation*}
	Then
	\begin{equation*}
		\supHdensity{m}{A \sim \mathbf{X}_{k,\alpha}(a,T,f,\kappa)}{a} < M(1-2^{-m})^{-1}
	\end{equation*}
	for every $ \kappa > 2^{k + \alpha}\lambda $.
\end{Lemma}

\begin{Proof}
	Arguing as in the proof of \ref{alternative characterization of approximate differentiability for functions} we conclude that 
	\begin{equation*}
		\supHdensity{m}{ A \cap \{ z : |f(T_{\natural}(z)) - T^{\perp}_{\natural}(z)| > 2^{k + \alpha}\lambda|z-a|^{k + \alpha} \}   }{a} < M(1-2^{-m})^{-1}.
	\end{equation*}
	Since $ \Der f(T_{\natural}(a)) = 0 $ we can easily get that 
	\begin{equation*}
		\lim_{r \to 0}\frac{ \Haus{m}( A \cap \mathbf{B}(a,r) \cap \{ z: |T^{\perp}_{\natural}(z-a)| > \epsilon\,r \})  }{\bm{\alpha}(m)\,r^{m}} = 0 \; \textrm{for every $ \epsilon > 0 $}
	\end{equation*}
	and applying \ref{remark on approximate tangent cone I} we conclude that $ \Hdensity{m}{A \sim \mathbf{X}(a,T,\epsilon)}{a} = 0 $. Since 
	\begin{flalign*}
		&  \mathbf{X}(a,T,\epsilon) \cap \{ z : |f(T_{\natural}(z)) - T^{\perp}_{\natural}(z)| \leq 2^{k + \alpha}\lambda|z-a|^{k + \alpha}\} \\
		& \quad \subseteq \mathbf{X}_{k, \alpha}(a,T,f,2^{k + \alpha}\lambda(1+\epsilon^{2})^{(k + \alpha)/2}) \quad \textrm{for every $ \epsilon > 0 $},
	\end{flalign*}
	the conclusion follows.
\end{Proof}

\begin{Definition}\label{ap diff for sets}
	Let $ n \geq 1 $ and $ k \geq 1 $ be integers, $ 0 \leq \alpha \leq 1 $, $ A \subseteq \Real{n} $, $ a \in \Real{n} $ and $ A_{1} = \{ x-a : x \in A \} $. We say that $ A $ is \textit{approximately differentiable of order $ (k,\alpha) $ at $ a $} if there exist an integer $ 1 \leq m \leq n $, \mbox{$ T \in \mathbf{G}(n,m) $} and a polynomial function $ P : T \rightarrow T^{\perp} $ of degree at most $ k $ such that $ P(0) = 0 $, $ \Der P(0) = 0 $ and the following two conditions hold.
	\begin{enumerate}
		\item \label{ap diff for sets: dimension} For every $ \epsilon> 0 $ there exists $ \rho > 0 $ and $ \eta > 0 $ such that 
		\begin{equation*}
			\Haus{m}(\mathbf{C}(T,z,\epsilon r , \epsilon r) \cap A_{1}) \geq \eta r^{m}
		\end{equation*}
		for every $ z \in T \cap \mathbf{B}(0,r) $ and $ 0 \leq r \leq \rho $.
		\item \label{ap diff for sets: approx by polynomials} For every $ \epsilon > 0 $
		\begin{equation*}
			\lim_{r \to 0} \frac{ \Haus{m}\left( A_{1} \cap \mathbf{B}(0,r) \cap \{ z: \bm{\delta}_{\gr(P)}(z) > \epsilon\,r^{k} \} \right) }{\bm{\alpha}(m)r^{m}  }  =0
		\end{equation*}
		and, if $ \alpha > 0 $, there exists $ 0 \leq \lambda < \infty $ such that
		\begin{equation*}
			\lim_{r \to 0} \frac{ \Haus{m}\left( A_{1} \cap \mathbf{B}(0,r) \cap \{ z: \bm{\delta}_{\gr(P)}(z) > \lambda \,r^{k+\alpha} \} \right) }{\bm{\alpha}(m)r^{m}  }  = 0.
		\end{equation*}
	\end{enumerate}
\end{Definition}

\begin{Remark}\label{Mattila tangent plane}
	If $ k =1 $ and $ \alpha = 0 $ the conditions in \ref{ap diff for sets} are equivalent to \cite[15.7]{MR1333890}.
\end{Remark}

\begin{Remark}\label{inclusion of approx tangent planes}
	We prove that the condition 
	\begin{equation*}
		T \subseteq \Tan^{m}_{*}(\Haus{m}\restrict A,a) 
	\end{equation*}
	is necessary and sufficient to have \ref{ap diff for sets}\eqref{ap diff for sets: dimension}. The condition is clearly necessary. To prove the sufficiency assume $ a = 0 $, suppose $ 0<\epsilon <1 $ and observe there exist an integer $ l \geq 1 $, $ v_{1}, \ldots , v_{l} \in \mathbf{S}^{n-1} \cap T $ and a positive number $ \eta $ such that 
	\begin{equation*}
		\textstyle T \cap \mathbf{S}^{n-1} \subseteq \bigcup_{i=1}^{l} \mathbf{U}(v_{i},\epsilon) \cap T,
	\end{equation*}
	\begin{equation*}
		\Haus{m}(A \cap \mathbf{U}(rv_{i},\epsilon r)) \geq \eta \epsilon^{-m} r^{m} \quad \textrm{whenever $ 0 < r \leq \eta $ and $ i = 1 , \ldots , l $.}
	\end{equation*}
	Since $ \infHdensity{m}{A}{0} > 0 $ by \ref{Remark on lower and upper densities and cones}, we can choose $ \eta > 0 $ smaller, if necessary, in order to have
	\begin{equation*}
		\Haus{m}(A \cap \mathbf{U}(0,\epsilon r)) \geq \eta r^{m} \quad \textrm{whenever $ 0 < r \leq \eta $.}
	\end{equation*}
	We fix $ 0 < r \leq \eta $ and $ z \in \mathbf{B}(0,r) $. If $ |z| \leq \epsilon r $ then $ \mathbf{U}(0,\epsilon r) \subseteq \mathbf{U}(z,2\epsilon r) $ and 
	\begin{equation*}
		\Haus{m}(A \cap \mathbf{U}(z,2 \epsilon r)) \geq \eta r^{m}.
	\end{equation*}
	If $ |z| \geq \epsilon r $ then we choose $ 1 \leq i \leq l $ such that $ |(z/|z|) - v_{i}| < \epsilon $ and we observe
	\begin{equation*}
		\mathbf{U}(|z|v_{i},|z|\epsilon) \subseteq \mathbf{U}(z, 2\epsilon |z|) \subseteq \mathbf{U}(z, 2\epsilon r),
	\end{equation*}
	\begin{equation*}
		\Haus{m}(A \cap \mathbf{U}(z,2 \epsilon r)) \geq \eta \epsilon^{-m} |z|^{m} \geq \eta r^{m}.
	\end{equation*}
\end{Remark}

\begin{Lemma}\label{distance vs vertical distance}
	Suppose $ 1 \leq m \leq n $ are integers, $ 0 < r < \infty $, $ w \in \Real{n} \cap \mathbf{B}(0,r) $, $ T \in \mathbf{G}(n,m) $ and $ f : T \rightarrow T^{\perp} $ is a locally Lipschitzian function such that $ f(0) = 0 $, . 
	
	Then $ \bm{\delta}_{\gr f}(w) \leq | T^{\perp}_{\natural}(w) - f (T_{\natural}(w))| \leq \big(2 + \Lip(f|\mathbf{B}(0,2r))\big)\bm{\delta}_{\gr f}(w) $.
\end{Lemma}

\begin{Proof}
	If we choose $ \chi \in T $ so that $ \bm{\delta}_{\gr f}(w) = |w - \chi - f(\chi)| $ then $ \chi \in \mathbf{B}(0,2r) $ and we get
	\begin{flalign*}
		\bm{\delta}_{\gr f}(w) & \leq | T^{\perp}_{\natural}(w) - f (T_{\natural}(w))| \\
		& \leq  |w - \chi - f(\chi)| + |\chi + f(\chi) - T_{\natural}(w) - f(T_{\natural}(w))|\\
		&   \leq \big(2 + \Lip(f|\mathbf{B}(0,2r))\big)\bm{\delta}_{\gr f}(w).
	\end{flalign*}
	
\end{Proof}

\begin{Lemma}\label{approximation by non homogeneous polynomial}
	Suppose $ 1 \leq m \leq n $ are integers, $ \gamma > 0 $, $ A \subseteq \Real{n} $, $ B \subseteq \Real{n} $ such that $ 0 \in \Clos B $, $ f : \Real{n} \rightarrow \Real{n} $ is an univalent map onto $ \Real{n} $ such that $ f(0) = 0 $ and $ f $ and $ f^{-1} $ are locally Lipschitzian maps.
	
	Then the following two conditions are equivalent.
	\begin{enumerate}
		\item \label{approximation by non homogeneous polynomial:1} For every $ \epsilon > 0 $ [for some $ 0 \leq \epsilon < \infty $] 
		\begin{equation*}
			\lim_{r \to 0} \frac{ \Haus{m}( A \cap \mathbf{B}(0,r) \cap \{ z: \bm{\delta}_{B}(z) > \epsilon\,r^{\gamma} \} ) }{\bm{\alpha}(m)r^{m}  }  =0.
		\end{equation*}
		\item \label{approximation by non homogeneous polynomial:2} For every $ \epsilon > 0 $ [for some $ 0 \leq \epsilon < \infty $]
		\begin{equation*}
			\lim_{r \to 0} \frac{\Haus{m}( f[A] \cap \mathbf{B}(0,r) \cap \{ z: \bm{\delta}_{f[B]}(z) > \epsilon r^{\gamma} \} )}   {\bm{\alpha}(m)r^{m}  }  =0.
		\end{equation*}
	\end{enumerate}
\end{Lemma}
\begin{Proof}
	Since $ f[\Clos B] = \Clos f[B] $ we assume $ B $ to be closed. Moreover if we prove one implication we immediately get the other one. Therefore we prove that \eqref{approximation by non homogeneous polynomial:1} implies \eqref{approximation by non homogeneous polynomial:2}. Suppose $ 1 < \Gamma < \infty $ is such that 
	\begin{eqnarray*}
		& \Gamma^{-1}|z-w| \leq |f(z)-f(w)| \leq \Gamma |z-w|,  &\\
		& \Gamma^{-1}|z-w| \leq |f^{-1}(z)-f^{-1}(w)| \leq \Gamma |z-w| &
	\end{eqnarray*}
	whenever $ z,w \in \mathbf{B}(0,2) $. Evidently it is enough to show that
	\begin{flalign*}
		& f^{-1}[ f[A]\cap \mathbf{B}(0,r/\Gamma^{2}) \cap \{ w : \bm{\delta}_{f[B]}(w) > \epsilon r^{\gamma} \}] \\
		& \quad \subseteq A \cap \mathbf{B}(0, r/\Gamma) \cap \{ w : \bm{\delta}_{B}(w) > \Gamma^{-1}\epsilon r^{\gamma} \}
	\end{flalign*}
	for $ \epsilon > 0 $ and $ 0 < r \leq 1 $. Suppose $ z \in f[A] \cap \mathbf{B}(0,r/\Gamma^{2}) $ such that $ \bm{\delta}_{f[B]}(z) > \epsilon r^{\gamma} $. Let $ w \in f[B] $ such that $ |f^{-1}(z) - f^{-1}(w) | = \bm{\delta}_{B}(f^{-1}(z)) $ and observe 
	\begin{eqnarray*}
		& \bm{\delta}_{B}(f^{-1}(z)) \leq |f^{-1}(z)|, \quad  |f^{-1}(w)| \leq 2|f^{-1}(z)| \leq 2 \Gamma |z| \leq 2 \Gamma^{-1}r \leq 2, &\\
		& |w| \leq \Gamma |f^{-1}(w)| \leq 2 r \leq 2, &\\
		& \bm{\delta}_{B}(f^{-1}(z)) \geq \Gamma^{-1}|z-w| \geq \Gamma^{-1}\bm{\delta}_{f[B]}(z) > \Gamma^{-1}\epsilon r^{\gamma}. &
	\end{eqnarray*}
\end{Proof}

\begin{Lemma}\label{uniqueness}
	Let $ 1 \leq m \leq n $ and $ k \geq 1 $ be integers, $ T \in \mathbf{G}(n,m) $ and let $ P : T \rightarrow T^{\perp} $ and $ Q : T \rightarrow T^{\perp} $ be polynomial functions of degree at most $ k $ such that $ P(0)=0 $ and $ \Der^{i} Q(0) = 0 $ for $ i = 0, \ldots , k-1 $. Suppose for every $ \epsilon > 0 $ there exists $ \rho > 0 $ such that 
	\begin{equation*}
		\mathbf{C}(T,z,\epsilon r,\epsilon r^{k}) \cap \{ w : \bm{\delta}_{\gr(Q)}(w) \leq \epsilon r^{k} \} \neq \varnothing
	\end{equation*}
	whenever $ z \in \gr(P) \cap \mathbf{B}(0,r) $ and $ 0 < r \leq \rho $. 
	
	Then $ P = Q $.
\end{Lemma}

\begin{Proof}
	Let $ 0 \leq c < \infty $ such that $ |P(\chi)| \leq c |\chi| $ whenever $ \chi \in T \cap \mathbf{B}(0,1) $. If $ 0 < \epsilon \leq 1 $ and $ 0 < \rho \leq 1 $ are as in the hypothesis, \mbox{$ \chi \in \mathbf{B}(0,(1+c)^{-1}\rho) \cap T $} and \mbox{$ z = \chi + P(\chi) $} then $ |z| \leq (1+c)|\chi| \leq \rho $. Therefore there exists
	\begin{equation*}
		w \in \mathbf{C}(T,z,\epsilon(1+c)|\chi|, \epsilon (1+c)^{k} |\chi|^{k})
	\end{equation*}
	such that $  \bm{\delta}_{\gr(Q)}(w) \leq \epsilon (1+c)^{k} |\chi|^{k} $. If $ y \in \gr(Q) $ is such that $ |w-y| = \bm{\delta}_{\gr(Q)}(w) $ then
	\begin{equation*}
		|T_{\natural}(y) - \chi| \leq | T_{\natural}(y-w)| + | T_{\natural}(w) - \chi | \leq 2\epsilon (1+c)^{k} |\chi|,
	\end{equation*}
	\begin{equation*}
		| P(\chi) - Q(\chi)| \leq |T^{\perp}_{\natural}(z) - T_{\natural}^{\perp}(w)| + | T^{\perp}_{\natural}(w) - T^{\perp}_{\natural}(y)| + |T^{\perp}_{\natural}(y) - Q(\chi)|,
	\end{equation*}
	and the Taylor's formula (see \cite[p.\ 46]{MR0257325}) implies
	\begin{equation*}
		Q(T_{\natural}(y)) - Q(\chi) = \sum_{i=1}^{k}  \langle (T_{\natural}(y) - \chi)^{i}/i! \odot \chi^{k-i}/(k-i)!, \Der^{k}Q(0) \rangle,
	\end{equation*}
	\begin{equation*}
		|Q(T_{\natural}(y)) - Q(\chi)|  \leq c_{1}\epsilon |\chi|^{k},
	\end{equation*}
	where $ c_{1} = \|\Der^{k}Q(0) \| \sum_{i=1}^{k}2^{i}(1+c)^{ki}/(i!(k-i)!) $. Therefore $ |Q(\chi) - P(\chi) | \leq (2(1+c)^{k}+c_{1})\epsilon |\chi|^{k} $ and the conclusion follows.
	
\end{Proof}

\begin{Theorem}\label{basic characterization}
	Suppose $ 1 \leq m \leq n $ and $ k \geq 1 $ are integers, $ 0 \leq \alpha \leq 1 $, $ A \subseteq \Real{n} $, $ a \in \Real{n} $, $ A_{1} = \{ x-a : x \in A\} $, $ T \in \mathbf{G}(n,m) $ and $ P : T \rightarrow T^{\perp} $ is a polynomial function of degree at most $ k $ such that $ P(0) = 0 $, $ \Der P(0) = 0 $.
	
	Then the following two conditions are equivalent.
	\begin{enumerate}
		\item \label{basic characterization:ap diff} $ T $ and $ P $ satisfy \ref{ap diff for sets}\eqref{ap diff for sets: dimension} and \ref{ap diff for sets}\eqref{ap diff for sets: approx by polynomials}.
		\item \label{basic characterization:ap diff by hom approx}  If $ P_{i}(\chi) = \langle \chi^{i}/i!,\Der^{i}P(0) \rangle $ for $ \chi \in T $ and $ i = 1 , \ldots , k $ and 
		\begin{equation*}
			A_{i} = \{ x - P_{i-1}(T_{\natural}(x)) : x \in A_{i-1} \} \quad \textrm{for $ i = 2, \ldots, k $},
		\end{equation*}
		then the following two conditions hold:
		\begin{enumerate}
			\item \label{basic characterization:ap diff by hom approx:dim} for every $ i = 1, \ldots , k $ and for every $ \epsilon> 0 $ there exist $ \rho > 0 $ and $ \eta > 0 $ such that 
			\begin{equation*}
				\Haus{m}(\mathbf{C}(T,z,\epsilon r , \epsilon r^{i}) \cap A_{i}) \geq \eta \bm{\alpha}(m)r^{m}
			\end{equation*}
			for every $ z \in \gr(P_{i}) \cap \mathbf{B}(0,r) $ and $ 0 \leq r \leq \rho $,
			\item  \label{basic characterization:ap diff by hom approx:integer case} for every $ i = 1 , \ldots , k $ and for every $ \epsilon > 0 $
			\begin{equation*}
				\lim_{r \to 0} \frac{ \Haus{m}\left( A_{i} \cap \mathbf{B}(0,r) \cap \{ z: \bm{\delta}_{\gr(P_{i})}(z) > \epsilon\,r^{i} \} \right) }{\bm{\alpha}(m)r^{m}  }  =0
			\end{equation*}
			and, if $ \alpha > 0 $, there exists $ 0 \leq \lambda < \infty $ such that
			\begin{equation*}
				\lim_{r \to 0} \frac{ \Haus{m}\left( A_{k} \cap \mathbf{B}(0,r) \cap \{ z: \bm{\delta}_{\gr(P_{k})}(z) > \lambda \,r^{k+\alpha} \} \right) }{\bm{\alpha}(m)r^{m}  }  = 0.
			\end{equation*}
		\end{enumerate}
		
	\end{enumerate}
	In this case $ P $ is uniquely determined by $ a $, $ A $ and $ k $,
	\begin{equation*}
		\infHdensity{m}{A}{a} > 0, \; \Tan^{m}_{*}(\Haus{m}\restrict A,a) = \Tan^{*m}(\Haus{m}\restrict A,a) = T
	\end{equation*}
	and $ A $ is approximately differentiable of order $ (l,\beta) $ whenever either $ l < k $ and $ 0 \leq \beta \leq 1 $ or $ l = k $ and $ 0 \leq \beta \leq \alpha $.
\end{Theorem}

\begin{Proof}
	Assume $ a = 0 $ and suppose $ \sup \big\{ 1, \sum_{j=1}^{k} 2^{j} \| \Der^{j} P(0) \|  \big\} < \Gamma < \infty $.
	
	For $ i = 1 , \ldots , k $ we define $ Q_{i} = \sum_{j=1}^{i}P_{j} $ and $ f_{i} : \Real{n} \rightarrow \Real{n} $ by
	\begin{equation*}
		f_{i}(x) = x - Q_{i}(T_{\natural}(x)) + P_{i}(T_{\natural}(x))  \quad \textrm{for $ x \in \Real{n} $}.
	\end{equation*}
	We observe that for every $ i = 1 , \ldots , k $ the map $ f_{i} $ is a diffeomorphism of \mbox{class $ \infty $} onto $ \Real{n} $ and, by induction over $ i $, one may easily prove that
	\begin{equation*}
		f_{i}[A]= A_{i}, \quad f_{i}[\gr Q_{i}] = \gr P_{i}.
	\end{equation*}
	
	Now, using \ref{approximation by non homogeneous polynomial}, it is easy to see that \eqref{basic characterization:ap diff by hom approx} implies \eqref{basic characterization:ap diff}. Henceforth, we \mbox{assume \eqref{basic characterization:ap diff}.} First we prove that 
	\begin{equation*}
		\bm{\delta}_{\gr P}(z) \geq \bm{\delta}_{\gr Q_{i}}(z) - \Gamma |z|^{i+1} \quad \textrm{for every $ i = 1 , \ldots , k-1 $ and $ z \in \mathbf{B}(0,1) $.}
	\end{equation*}
	In fact if $ w \in \gr(P) $ such that $ |z-w| = \bm{\delta}_{\gr(P)}(z) $ then 
	\begin{flalign*}
		& |w| \leq 2 |z|, \quad |P(T_{\natural}(w)) - Q_{i}(T_{\natural}(w))| \leq \Gamma |z|^{i+1},\\
		|z-w| & \geq |z-T_{\natural}(w) - Q_{i}(T_{\natural}(w))| - |Q_{i}(T_{\natural}(w)) - P(T_{\natural}(w))| \\
		& \geq \bm{\delta}_{\gr Q_{i}}(z) - \Gamma |z|^{i+1}.
	\end{flalign*}
	It follows, for $ i = 1 , \ldots , k $, that
	\begin{equation*}
		\lim_{r \to 0} \frac{ \Haus{m}\left( A_{1} \cap \mathbf{B}(0,r) \cap \{ z: \bm{\delta}_{\gr(Q_{i})}(z) > \epsilon\,r^{i} \} \right) }{\bm{\alpha}(m)r^{m}  }  =0 \quad \textrm{for every $ \epsilon > 0 $.}
	\end{equation*}
and, for $ i = 1, \ldots, k-1 $ that
\begin{equation*}
 \lim_{r \to 0} \frac{ \Haus{m}\left( A_{1} \cap \mathbf{B}(0,r) \cap \{ z: \bm{\delta}_{\gr(Q_{i})}(z) > 2 \Gamma r^{i+1} \} \right) }{\bm{\alpha}(m)r^{m}  }  = 0.
\end{equation*}
Therefore \eqref{basic characterization:ap diff by hom approx:integer case} is a consequence of \ref{approximation by non homogeneous polynomial}. Moreover, using \ref{approximation by non homogeneous polynomial}, it follows that $ A $ is approximately differentiable of order $ (l,\beta) $ whenever either $ l < k $ and $ 0 \leq \beta \leq 1 $ or $ l = k $ and $ 0 \leq \beta \leq \alpha $. We prove now \eqref{basic characterization:ap diff by hom approx:dim}, whose proof is slightly more involved. We fix $ 2 \leq i \leq k $ and we replace $ \Gamma $ by a larger number, if necessary, in order to have
	\begin{equation*} 
		|f_{i}^{-1}(w) - f_{i}^{-1}(z)| \leq \Gamma |w-z|, \quad |f_{i}(w) - f_{i}(z)| \leq \Gamma |w-z| 
	\end{equation*}
	for $ w,z \in \mathbf{B}(0,4) $. Therefore we have
	\begin{equation*}
		f_{i}[A \cap \mathbf{C}(T,z,\epsilon r / \Gamma, \epsilon r/ \Gamma)] \subseteq A_{i} \cap \mathbf{C}(T,z,\epsilon r/\Gamma, 3r),
	\end{equation*}
	\begin{equation*}
		\Haus{m}(A_{i}\cap \mathbf{C}(T,z,\epsilon r/\Gamma, 3r)) \geq \Gamma^{-m}\Haus{m}(A \cap \mathbf{C}(T,z,\epsilon r/\Gamma, \epsilon r/\Gamma)),
	\end{equation*}
	whenever $ 0 < r \leq 1 $, $ 0 < \epsilon \leq 1 $ and $ z \in \mathbf{B}(0,r/\Gamma) $. We fix $ 0 < \epsilon \leq 1 $ and, using \ref{distance vs vertical distance}, we can find $ 0 < \rho \leq 1 $ and $ \eta > 0 $ such that 
	\begin{equation*}
		\Haus{m}(A_{i} \cap \mathbf{B}(0,5r) \cap \{w : |T_{\natural}^{\perp}(w) - P_{i}(T_{\natural}(w))| > \epsilon r^{i} \}) < \eta \Gamma^{-m} \bm{\alpha}(m) r^{m},
	\end{equation*}
	\begin{equation*}
		\Haus{m}(\mathbf{C}(T,z,\epsilon r/\Gamma , \epsilon r/\Gamma) \cap A) \geq 2 \eta \bm{\alpha}(m) r^{m} \quad \textrm{for every $ z \in T \cap \mathbf{B}(0,r) $},
	\end{equation*}
	whenever $ 0 < r \leq \rho $. Let $ 0 < r \leq \rho $ and $ z \in \gr(P_{i}) \cap \mathbf{B}(0,r/\Gamma) $. Then
	\begin{flalign*}
		\mathbf{C}(T,z,\epsilon r / \Gamma,3r) & \subseteq \big(\mathbf{B}(0,5r) \cap \{ w : | P_{i}(T_{\natural}(w)) - T^{\perp}_{\natural}(w)| > \epsilon r^{i} \}\big)\\
		& \quad \cup \mathbf{C}(T,z,\epsilon r / \Gamma,2\epsilon r^{i}).
	\end{flalign*}
	In fact if $ w \in \mathbf{C}(T,z,\epsilon r/\Gamma , 3r ) $ and $ | P_{i}(T_{\natural}(w)) - T^{\perp}_{\natural}(w)| \leq \epsilon r^{i} $ then
	\begin{flalign*}
		& | P_{i}(T_{\natural}(w)) - P_{i}(T_{\natural}(z))|  \\
		& \quad =\textstyle \big|\sum_{j=1}^{i} \langle (T_{\natural}(w-z)^{j}/j!) \odot (T_{\natural}(z)^{i-j}/(i-j)!) , \Der^{i}P(0) \rangle \big| \\
		& \quad \leq i \| \Der^{i}P(0) \| \Gamma^{-i} \epsilon r^{i}  \leq \epsilon r^{i}
	\end{flalign*}
	and we infer
	\begin{equation*}
		|T_{\natural}^{\perp}(z) - T_{\natural}^{\perp}(w)| \leq 2 \epsilon r^{i}, \quad w \in \mathbf{C}(T,z,\epsilon r/\Gamma, 2\epsilon r^{i}).
	\end{equation*}
	We can now conclude that
	\begin{equation*}
		\Haus{m}(A_{i}\cap \mathbf{C}(T,z,\epsilon r/\Gamma, 2\epsilon r^{i})) \geq \Gamma^{-m}\eta \bm{\alpha}(m)r^{m}
	\end{equation*}
	and \eqref{basic characterization:ap diff by hom approx:dim} is proved.
	
	By \ref{ap diff for sets}\eqref{ap diff for sets: dimension} we immediately conclude that $ \infHdensity{m}{A}{0} > 0 $ and $ T \subseteq \Tan^{m}_{*}(\Haus{m}\restrict A,0)  $. By \eqref{basic characterization:ap diff by hom approx:integer case} and \ref{remark on approximate tangent cone I} we conclude that $ \Tan^{*m}(\Haus{m}\restrict A,0) \subseteq T $. Finally let $ R : T \rightarrow T^{\perp} $ be a polynomial function of degree at most $ k $ such that $ R(0) = 0 $ and $ \Der R(0) = 0 $ and satisfying \ref{ap diff for sets}\eqref{ap diff for sets: dimension} and \ref{ap diff for sets}\eqref{ap diff for sets: approx by polynomials}. Let 
	\begin{eqnarray*}
		& R_{i}(\chi) = \langle \chi^{i}/i!,\Der^{i}R(0) \rangle \quad \textrm{for $ \chi \in T $ and $ i = 1 , \ldots , k $}, &\\
		& B_{1} = A_{1}, \quad B_{i} = \{ x - R_{i-1}(T_{\natural}(x)) : x \in B_{i-1} \} \; \textrm{for $ i = 2, \ldots, k $}. &
	\end{eqnarray*}
	We prove by induction that $ P_{i} = R_{i} $ for $ i = 1 , \ldots , k $. Assume, for $ j = 1 , \ldots , i $ and $ i < k $, that $ P_{j} = R_{j} $  and observe that $ A_{i+1} = B_{i+1} $. Let $ \epsilon > 0 $, $ 0 < \rho \leq 1 $ and $ \eta > 0 $ such that
	\begin{equation*}
		\Haus{m}(\mathbf{C}(T,z,\epsilon r , \epsilon r^{i+1}) \cap B_{i+1}) \geq \eta \bm{\alpha}(m)r^{m} \; \textrm{for every $ z \in \mathbf{B}(0,r) \cap \gr(R_{i+1}) $},
	\end{equation*}
	\begin{equation*}
		\Haus{m}\left( A_{i+1} \cap \mathbf{B}(0,2r) \cap \{ z: \bm{\delta}_{\gr(P_{i+1})}(z) > \epsilon\,r^{i+1} \} \right) \leq (\eta/2)\bm{\alpha}(m) r^{m}
	\end{equation*}
	whenever $ 0 < r \leq \rho $. Therefore for every $ z \in \mathbf{B}(0,r) \cap \gr(R_{i+1}) $ and for every $ 0 < r \leq \rho $ we conclude that 
	\begin{equation*}
		\Haus{m}(B_{i+1} \cap \mathbf{C}(T,z,\epsilon r , \epsilon r^{i+1}) \cap \{ z: \bm{\delta}_{\gr(P_{i+1})}(z) \leq \epsilon\,r^{i+1} \}  ) \geq (\eta/2) \bm{\alpha}(m)r^{m}
	\end{equation*}
	and $ P_{i+1} = R_{i+1} $ by \ref{uniqueness}.
\end{Proof}

\begin{Remark}
	A conceptually similar characterization has been proved for the notion of pointwise differentiability in \cite[3.22]{snulmenn:sets.v1}. Moreover the reader may find useful to compare \ref{basic characterization}\eqref{basic characterization:ap diff by hom approx} and \cite[3.4]{MR1285779}, where a concept of approximate tangent paraboloid is introduced by means of inhomogeneous dilations and weak convergence of Radon measures.
\end{Remark}

\begin{Remark}\label{alternative characterization of approximately differentiability}
	Suppose $ A \subseteq \Real{n} $ and $ a \in \Real{n} $. It is not difficult to see that the condition
	\begin{equation*}
		\Tan^{m}(\Haus{m}\restrict A,a)  \in \mathbf{G}(n,m) \; \textrm{for some integer $ 1 \leq m \leq n $}
	\end{equation*}
	is necessary and sufficient to conclude that $ A $ is approximately differentiable of order $ 1 $ at $ a $. In fact the necessity is asserted in \ref{basic characterization}, while the sufficiency follows from \ref{inclusion of approx tangent planes} and \ref{remark on approximate tangent cone I}.
\end{Remark}

\begin{Remark}\label{approximate differentiable set with infinite density}
	We describe now a simple example which illustrates some features of the notion of approximate differentiability of order $ 1 $. 
	
	With each $ \gamma > 1 $ and $ \gamma^{-1} < \alpha < (\gamma-1)^{-1} $ we associate the family $ F_{\alpha, \gamma} $ consisting of the subsets
	\begin{equation*}
		\Real{2} \cap \big( \{ (n^{-\alpha}, t) : 0 \leq t \leq n^{-\alpha\gamma} \} \cup \{ (-n^{-\alpha}, t) : 0 \leq t \leq n^{-\alpha\gamma} \}\big)
	\end{equation*}
	correspoding to the integers $ n \geq 1 $. We define
	\begin{equation*}
		A_{\alpha, \gamma} =   \big(\Real{2}\cap \{ (s,0) : -1 \leq s \leq 1 \} \big) \cup \bigcup F_{\alpha, \gamma}.
	\end{equation*}
	Since $ \alpha \gamma > 1 $ then $ \Haus{1}(A_{\alpha, \gamma}) < \infty $. Moreover, for each $ n \geq 1 $, 
	\begin{equation*}
		(n-1)^{\alpha} \sum_{i=n}^{\infty} i^{-\alpha\gamma} \geq (n-1)^{\alpha}\int_{n}^{\infty}x^{-\alpha\gamma}d\Leb{1}x = (n-1)^{\alpha}(\alpha\gamma-1)^{-1} n^{1-\alpha\gamma} \to \infty 
	\end{equation*}
	as $ n \to \infty $. Therefore $ \infHdensity{1}{A_{\alpha,\gamma}}{0} = \infty $. Finally $ A_{\alpha, \gamma} $ is approximately differentiable of order $ 1 $ at $ 0 $ by \ref{alternative characterization of approximately differentiability}, since 
	\begin{equation*}
		\Tan^{1}(\Haus{1}\restrict A_{\alpha, \gamma},0) = \Real{} \times \{0\}.
	\end{equation*}
	
\end{Remark}

\begin{Remark}\label{uniqueness of the dimension integer}
	Let $ A \subseteq \Real{n} $, $ a \in \Real{n} $ and let $ 0 \leq \mu \leq \nu $ be integers. Since 
	\begin{equation*}
		\Tan^{*\nu}(\Haus{\nu}\restrict A,a) \subseteq \Tan^{*\mu}(\Haus{\mu}\restrict A,a), 
	\end{equation*}
	we deduce by \ref{basic characterization} that the integer $ m $ in \ref{ap diff for sets} is uniquely determined by $ A $ \mbox{and $ a $.} 
\end{Remark}

\begin{Definition}\label{definition of approximate tangent space}
	Let $ A \subseteq \Real{n} $ and $ a \in \Real{n} $. Suppose $ A $ is approximately differentiable of order $ 1 $ at $ a $ and $ m $ and $ T $ are as \mbox{in \ref{ap diff for sets}.} We define \textit{the approximate tangent space of $ A $ at $ a $} to be the $ m $ dimensional subspace $ T $ and we denote it by $ \ap \Tan(A,a) $. Moreover we define \textit{the approximate normal space of $ A $ at $ a $} to be 
	\begin{equation*}
		\ap \Nor(A,a) = \Real{n} \cap \{ v : v \bullet u = 0 \}.
	\end{equation*}
	
\end{Definition}

\begin{Definition}\label{definition of approximate differentials}
	Let $ A \subseteq \Real{n} $, let $ k \geq 2 $ be an integer and $ a \in \Real{n} $. If $ A $ is approximately differentiable of \mbox{order $ k $} at $ a $ then we define \textit{the approximate differential of order $ k $ of $ A $ at $ a $} to be the symmetric $ k $ linear map 
	\begin{equation*}
		\ap \Der^{k}A(a) = \Der^{k}(P\circ T_{\natural})(0)\in \textstyle \bigodot^{k}(\Real{n},\Real{n}),
	\end{equation*}
	where $ T = \ap \Tan(A,a) $ and $ P : T \rightarrow T^{\perp} $ is as in \ref{ap diff for sets}.
\end{Definition}

\begin{Remark}\label{comparison with Simon approximate tangent plane}
	Suppose $ A \subseteq \Real{n} $.
	
	Following \cite[11.2, 11.4]{MR756417} (see also \cite[2.2]{MR1686704}) we consider the map $ P_{A} $ whose domain is given by the set of $ a \in \Real{n} $ such that there exist an integer $ 1 \leq m \leq n $, $ T \in \mathbf{G}(n,m) $ and $ 0 < \theta < \infty $ such that
	\begin{equation*}
		\lim_{r \to 0+}r^{-m} \int_{A}f((x-a)/r)d\Haus{m}x = \theta \int_{T}f d\Haus{m} \quad \textrm{whenever $ f \in \mathscr{K}(\Real{n}) $},
	\end{equation*}
	and whose value $ P_{A}(a) $ at $ a $ equals $ T $. In fact, one may readily verifies that $ m $, $ T $ \mbox{and $ \theta $} are uniquely determined by $ A $ and $ a $.
	
	Then it is not difficult to check that if $ a \in \dmn P_{A} $ and $ m = \dim P_{A}(a) $ then 
	\begin{equation*}
		P_{A}(a) \subseteq \Tan_{*}^{m}(\Haus{m} \restrict A,a), \quad  \Tan^{*m}(\Haus{m} \restrict A,a) \subseteq P_{A}(a),
	\end{equation*}
	\begin{equation*}
		\Hdensity{m}{A}{a} = \theta.
	\end{equation*}
	
	Using \ref{alternative characterization of approximately differentiability} and \ref{basic characterization} we deduce
	\begin{equation*}
		\dmn P_{A} \subseteq \dmn \ap \Tan(A,\cdot), \quad P_{A}(a) = \ap \Tan(A,a) \; \textrm{whenever $ a \in \dmn P_{A} $}.
	\end{equation*}
	If $ A_{\alpha, \gamma} $ is defined as in \ref{approximate differentiable set with infinite density}, then $ 0 \in (\dmn \ap \Tan(A_{\alpha,\gamma},\cdot)) \sim (\dmn P_{A_{\alpha,\gamma}}) $.
\end{Remark}

\begin{Remark}\label{symmetric difference with ap diff sets}
	Let $ 1 \leq m \leq n $ and $ k \geq 1 $ be integers, $ 0 \leq \alpha \leq 1 $, $ A \subseteq \Real{n} $, $ B \subseteq \Real{n} $ and $ a \in \Real{n} $. Suppose $ A $ is approximately differentiable of order $ (k,\alpha) $ at $ a $, $ m = \dim \ap \Tan(A,a) $ and 
	\begin{equation*}
		\Hdensity{m}{A \sim B}{a}= 0, \quad \Hdensity{m}{B \sim A}{a}= 0.
	\end{equation*}
	Then $ B $ is approximately differentiable of order $ (k,\alpha) $ at $ a $ with 
	\begin{eqnarray*}
		& \ap \Tan(A,a) = \ap \Tan(B,a), &\\
		& \ap\Der^{i}A(a) = \ap\Der^{i}B(a) \quad \textrm{for $ i = 2, \ldots,k $.} &
	\end{eqnarray*}
\end{Remark}

\begin{Theorem}\label{approximate differentiability of rectifiable sets}
	Let $ 1 \leq m \leq n $ and $ k \geq 1 $ be integers, $ 0 \leq \alpha \leq 1 $ and let $ A \subseteq \Real{n} $ be $ \Haus{m} $ measurable and $ \rect{m} $ rectifiable of class $ (k,\alpha) $. 
	
	Then for $ \Haus{m} $ a.e.\ $ a \in A $ the set $ A $ is approximately differentiable of order $ (k,\alpha) $ at $ a $ with
	\begin{equation*}
		\ap \Tan(A,a) \in \mathbf{G}(n,m) 
	\end{equation*}
\end{Theorem}

\begin{Proof}
	Since an $ m $ dimensional submanifold $ M $ of class $ (k,\alpha) $ of $ \Real{n} $ locally corresponds at each $ a \in M  $ to a graph of function $ f : \Tan(M,a) \rightarrow \Nor(M,a) $ of class $ (k,\alpha) $ with $ \Der f(\Tan(M,a)_{\natural}(a)) = 0 $, one readily checks that $ M $ is approximately differentiable of order $ (k,\alpha) $ at each of its points. Then the conclusion follows from \cite[2.10.19(4)]{MR0257325} and \ref{symmetric difference with ap diff sets}.
\end{Proof}

\begin{Remark}
	The conclusion of \ref{approximate differentiability of rectifiable sets} may not hold if we replace ``$ \rect{m} $ rectifiable'' with ``countably $ \rect{m} $ rectifiable'', as the following example in $ \Real{2} $ shows,
	\begin{equation*}
	\textstyle	\bigcup_{n=1}^{\infty} \{ (n^{-1},t) : 0 \leq t \leq 1 \} \cup \{ (0,t): 0 \leq t \leq 1  \}.
	\end{equation*}
\end{Remark}

\begin{Theorem}\label{rectifiable sets and normal vector fields}
	
	Let $ 1 \leq m \leq n $ be integers, let $ A \subseteq \Real{n} $ be $ \Haus{m} $ measurable and $ \rect{m} $ rectifiable of class $ 1 $ and let $ \nu: A \rightarrow \Real{n} $ be a map such that for $ \Haus{m} $ a.e.\ $ x \in A $ there exists $ 0 \leq \lambda < \infty $ such that \
	\begin{equation*}
		\Hdensity{m}{A \cap \{ z: |\nu(z)-\nu(x)| > \lambda |z-x| \}}{x} = 0.
	\end{equation*}
	
	Then $ \nu $ is $ \Haus{m}\restrict A $ measurable and $ (\Haus{m}\restrict A,m) $ approximately differentiable\footnote{See \cite[3.2.16]{MR0257325}.} at $ \Haus{m} $ a.e.\ $ x \in A $.
	
	If additionally $ \nu(x) \in \ap \Nor(A,x) $ for $ \Haus{m} $ a.e.\ $ x \in A $ and $ A $ is $ \rect{m} $ rectifiable of class $ 2 $ then
	\begin{equation*}
		(\Haus{m}\restrict A,m)\ap\Der \nu(x)(u) \bullet v = - \ap \Der^{2} A(x)(u,v) \bullet \nu(x)
	\end{equation*}
	for every $ u,v \in \ap \Tan(A,x) $ and for $ \Haus{m} $ a.e.\ $ x \in A $.
\end{Theorem}

\begin{Proof}
	By \cite[3.2.29, 3.1.19(4), 2.10.19(4), 3.2.16]{MR0257325} it is enough to prove the statement in the following special case: let $ U \subseteq \Real{n} $, $ V \subseteq \Real{m} $ be bounded open sets and let $ \phi : U \rightarrow \Real{m} $, $ \psi : V \rightarrow \Real{n} $ be maps of class $ 1 $ (of class $ 2 $ if $ A $ is $ \rect{m} $ rectifiable of class $ 2 $) such that $ A \subseteq \im \psi $ and $ \phi \circ \psi = \mathbf{1}_{V} $. Let $ M = \im \psi $ and observe that $ \phi|M = \psi^{-1} $, $ \phi[A] $ is an \mbox{$ \Haus{m} $ measurable} subset of $ \Real{m} $. Moreover we can prove that $ \nu $ is $ \Haus{m}\restrict A $ measurable by \cite[2.9.13]{MR0257325}. In fact one verifies that $ V = \{ (a,\mathbf{B}(a,r)): a \in \Real{n}, \; 0 < r < \infty \} $ is an $ \Haus{m}\restrict A $ Vitali relation by \cite[2.8.18]{MR0257325} and, since $ \Hdensity{m}{A}{x} = 1 $ for $ \Haus{m} $ a.e. $ x \in A $ by \cite[3.2.19]{MR0257325}, we conclude that $ \nu $ is $ (\Haus{m}\restrict A,V) $ approximately continuous\footnote{See \cite[2.9.12]{MR0257325}.} at $ \Haus{m} $ a.e.\ $ x \in A $.
	
	Let $ \eta = \nu \circ (\psi| \phi[A]) $ and, by \cite[2.9.11]{MR0257325}, we deduce that $ \eta $ is approximately differentiable of order $ (0,1) $ at $ \Leb{m} $ a.e. $ \chi \in \phi[A] $. Therefore by \ref{rectifiability and approximate differentiability for functions I}\eqref{rectifiability and approximate differentiability for functions I: rectifiability} there exist countably many maps $ \eta_{j} : \Real{m} \rightarrow \Real{n} $ of class $ 1 $ such that 
	\begin{equation*}
		\textstyle \Leb{m}\left( \phi[A] \sim \bigcup_{j=1}^{\infty}\{ \chi : \eta_{j}(\chi) = \eta(\chi) \} \right) = 0.
	\end{equation*}
	We deduce, by \cite[2.10.19(4)]{MR0257325}, that $ \nu $ is $ (\Haus{m}\restrict A,m) $ approximately differentiable at $ \Haus{m} $ a.e.\ $ x \in A $ because 
	\begin{equation*}
		\textstyle \Haus{m}\left(A \sim \bigcup_{j=1}^{\infty}\{ x : (\eta_{j} \circ \phi)(x) = \nu(x) \} \right) = 0.
	\end{equation*}
	
	If we further assume $ \nu(x) \in \ap \Nor(A,x) $ for $ \Haus{m} $ a.e.\ $ x \in A $ and $ A $ is $ \rect{m} $ rectifiable of class $ 2 $ then, for every $ j \geq 1 $, we define
	\begin{equation*}
		\nu_{j}(x) = \big(\Nor(M,x)_{\natural} \circ \eta_{j} \circ \phi\big)(x) \quad \textrm{for $ x \in M $},
	\end{equation*}
	we observe that $ \nu_{j} $ is of class $ 1 $ relative to $ M $ and, by \cite[2.10.19(4)]{MR0257325} \mbox{and \ref{symmetric difference with ap diff sets},}
	\begin{equation*}
		\textstyle \Haus{m}\left(A \sim \bigcup_{j=1}^{\infty}\{ x : \nu_{j}(x) = \nu(x) \} \right) = 0.
	\end{equation*}
	Since, by \cite[2.10.19(4), 3.2.16]{MR0257325} and \ref{symmetric difference with ap diff sets},
	\begin{equation*}
		\Der \nu_{j}(x)(u) \bullet v = - \ap \Der^{2}A(x)(u,v) \bullet \nu_{j}(x) \quad \textrm{for every $ u,v \in \ap \Tan(A,x) $}
	\end{equation*}
	and $ (\Haus{m}\restrict A,m)\ap\Der \nu(x) = \Der \nu_{j}(x) $ for $ \Haus{m} $ a.e.\ $ x \in A $, the conclusion follows.
\end{Proof}

\begin{Remark}
	The conclusion of the second part of \ref{rectifiable sets and normal vector fields} may fail to hold at $ \Haus{m} $ a.e.\ $ a \in A $ if we omit the hypothesis ``$ A $ is $ \rect{m} $ rectifiable of class $ 2 $'', even if we assume that $ A $ is an $ m $ dimensional submanifold of class $ 1 $. This fact can be easily deduced from \cite{MR0427559} and \ref{the set of approximate differentiable points is rectifiable}.
	
	Moreover the same conclusion may fail to hold at $ \Haus{m} $ a.e.\ $ a \in A $ if we omit the hypothesis ``$ A $ is $ \rect{m} $ rectifiable of class $ 2 $'' but we assume $ \nu(x) = \zeta $ for $ \Haus{m} $ a.e. $ x \in A $ for some $ \zeta \in \mathbf{S}^{n-1} $. In fact it is proved in \cite[Appendix]{MR1285779} that for every $ 0 < \alpha < 1 $ there exists a function $ f : \Real{} \rightarrow \Real{} $ of class $ (1,\alpha) $ and a Cantor-type set $ E \subseteq \Real{} $ such that 
	\begin{eqnarray*}
		&	\Leb{1}(E) > 0, \quad \Der f (x) = 0 \quad \textrm{for every $ x \in E $}, &\\
		& \Leb{1}\big(E \cap \{x : f(x)= g(x)\} \big) = 0 \quad \textrm{whenever $ g : \Real{} \rightarrow \Real{} $ is of class $ 2 $}. &
	\end{eqnarray*}
	If $ A = \gr(f|E) $ then, by , $ \Haus{1}( A \cap \dmn \ap \Der^{2}A) = 0 $.
\end{Remark}

\section{Relation with pointwise differentiability}\label{section Relation with pointwise differentiability}

\begin{Definition}\label{def pf pt tangent cone}
	Suppose $ X $ is a normed vector space, $ B \subseteq X $ and $ a \in X $.
	
	We define \textit{the upper tangent cone of $ B $ at $ a $} by
	\begin{equation*}
		\Tan^{*}(B,a) = X \cap \{ v : \liminf_{r \to 0+} r^{-1} \bm{\delta}_{B}(a+rv) = 0 \}
	\end{equation*}
	and \textit{the lower tangent cone of $ B $ at $ a $} by
	\begin{equation*}
		\Tan_{*}(B,a) = X \cap \{ v : \lim_{r \to 0+} r^{-1} \bm{\delta}_{B}(a+rv) = 0 \}
	\end{equation*}
	In case $\Tan_{*}(B,a)= \Tan^{*}(B,a) $, this set is denoted by $ \Tan(B,a) $ and we call it \textit{the tangent cone of $ B $ at $ a $}. Finally \textit{the [lower, upper] normal cone of $ B $ at $ a $} is defined to be the set of $ v \in \Real{n} $ such that $ v \bullet u \leq 0 $ whenever \mbox{$ u \in [\Tan_{*}(A,a), \; \Tan^{*}(A,a)] \; \Tan(A, a) $}, and it is denoted by 
\begin{equation*}
[\Nor_{*}(A,a), \;  \Nor^{*}(A,a)] \;  \Nor(A,a).
\end{equation*}                             
                                            
\end{Definition}

\begin{Remark}\label{inclusion of tangent cones}
	If $ 1 \leq m \leq n $ are integers and $ B \subseteq \Real{n} $ then one may verify that
	\begin{eqnarray*}
		& \Tan^{*m}(\Haus{m}\restrict B,a) \subseteq \Tan^{*}(B,a) & \\
		& \quad \rsubseteq \qquad \qquad \qquad \quad \rsubseteq \quad \\
		& \Tan_{*}^{m}(\Haus{m}\restrict B,a)  \subseteq \Tan_{*}(B,a). &
	\end{eqnarray*}
	
	Moreover one may readily verify that $ \Tan_{*}(B,a) $ and $ \Tan^{*}(B,a) $ are closed cones.
\end{Remark}

\begin{Remark}
	This notation does not agree with \cite[4.3]{MR0110078}, \cite[3.1.21]{MR0257325} and \cite{snulmenn:sets.v1}. In fact $ \Tan^{*}(B,a) $ is denoted by $ \Tan(B,a) $ therein.
\end{Remark}

\begin{Remark}
 Employing \cite[4.1]{MR0110078}, we observe that if $ A $ is a closed subset \mbox{of $ \Real{n} $,} $ a \in A $ and $ \reach(A,a) > 0 $ then, by \cite[4.8(10),(12)]{MR0110078}, 
\begin{equation*}
 \Tan_{*}(A,a) = \Tan^{*}(A,a).
\end{equation*}
\end{Remark}

\begin{Definition}\label{definition of pointwise differentiability}
	Let $ k $ and $ n $ be positive integers, $ 0 \leq \alpha \leq 1 $ and $ B \subseteq \Real{n} $. We say that \textit{$ B $ is pointwise differentiable of order $ (k,\alpha) $ at $ a $} if there exists a submanifold $ M \subseteq \Real{n} $ of class $ (k,\alpha) $ such that $ a \in M $,
	\begin{equation*}
		\lim_{r \to 0} r^{-1} \sup\{ |\bm{\delta}_{M}(x) - \bm{\delta}_{B}(x)| : x \in \mathbf{B}(a,r) \} = 0,
	\end{equation*}
	\begin{equation*}
		\lim_{r \to 0} r^{-k} \sup \{ \bm{\delta}_{M}(x): x \in \mathbf{B}(a,r) \cap B \} = 0 \quad \textrm{if $ \alpha = 0 $},
	\end{equation*}
	\begin{equation*}
		\limsup_{r \to 0} r^{-k-\alpha} \sup \{ \bm{\delta}_{M}(x): x \in \mathbf{B}(a,r) \cap B \} < \infty \quad \textrm{if $ \alpha > 0 $}.
	\end{equation*}
\end{Definition}

\begin{Remark}
	This concept has been introduced in \cite[3.3]{snulmenn:sets.v1}. In \ref{ap diff and pt diff} and \ref{touching by balls and estimate ssf} we employ the concept of pointwise differential of order $ i $ for sets, introduced in \cite[3.12]{snulmenn:sets.v1}.
\end{Remark}

\begin{Remark}
	It is worth to mention that, for sets, pointwise differentiability does not imply approximate differentiability. In fact, suppose $ n \geq 1 $ is an integer and $ B $ is a countable dense subset of $ \Real{n} $. Then for every integer $ k \geq 1 $ the set $ B $ is pointwise differentiable of order $ k $ at every $ x \in \Real{n} $. But $ B $ is not approximately differentiable of order $ 1 $ at every $ x \in \Real{n} $.
\end{Remark}

\begin{Lemma}\label{basic remark on pt tangent cone}
	Let $ B \subseteq \Real{n} $ and $ a \in \Clos B $. 
	
	Then the following statements hold.
	\begin{enumerate}
		\item \label{basic remark on pt tangent cone:1} If $ M = \{ a + v : v \in \Tan^{*}(B,a) \} $ then
		\begin{equation*}
			\lim_{r \to 0} r^{-1} \sup \{ \bm{\delta}_{M}(x) : x \in \mathbf{B}(a,r) \cap B \} = 0.
		\end{equation*}
		\item \label{basic remark on pt tangent cone:2}If $ M = \{ a + v : v \in \Tan_{*}(B,a) \} $ then
		\begin{equation*}
			\lim_{r \to 0} r^{-1} \sup \{ \bm{\delta}_{B}(x) : x \in \mathbf{B}(a,r) \cap M \} = 0.
		\end{equation*}
		\item \label{basic remark on pt tangent cone:3} The condition 
		\begin{equation*}
			\Tan(B,a) \in \mathbf{G}(n,m) \quad \textrm{for some integer $ 0 \leq m \leq n $}
		\end{equation*}
		is necessary and sufficient to conclude that $ A $ is pointwise differentiable of order $ 1 $ at $ a $.
	\end{enumerate}
\end{Lemma}

\begin{Proof}
	\textit{Proof of \eqref{basic remark on pt tangent cone:1}}. If there existed $ \epsilon> 0 $, $ r_{i} > 0 $, $ r_{i} \to 0 $ as $ i \to \infty $ and $ x_{i} \in B \cap \mathbf{B}(a,r_{i}) $ such that $ \bm{\delta}_{M}(x_{i}) \geq \epsilon r_{i} $ then, possibly passing to a subsequence, we could assume there would exist $ v \in \mathbf{S}^{n-1} $ such that $ (x_{i} -a) / |x_{i} -a| \to v $ as $ i \to \infty $. Then $ v \in \Tan^{*}(B,a) $,
	\begin{equation*}
		\epsilon \leq r_{i}^{-1}\big|x_{i}-a-|x_{i}-a|v \big| \leq |x_{i}-a|^{-1}\big|x_{i}-a-|x_{i}-a|v\big| \;\; \textrm{for $ i \geq 1 $}
	\end{equation*}
	and we would get a contradiction.
	
	\textit{Proof of \eqref{basic remark on pt tangent cone:2}}. Suppose $ \epsilon > 0 $ and observe there exist an integer $ l \geq 1 $, $ v_{1}, \ldots , v_{l} \in \Tan_{*}(B,a) \cap \mathbf{S}^{n-1} $ and $ \eta > 0 $ such that $ r^{-1}\bm{\delta}_{B}(a+rv_{i}) < \epsilon $ whenever $ i = 1, \ldots , l $ and $ 0 < r \leq \eta $ and 
	\begin{equation*}
		\textstyle \Tan_{*}(B,a) \cap \mathbf{S}^{n-1} \subseteq \bigcup_{i=1}^{l}\mathbf{B}(v_{i}, \epsilon).
	\end{equation*}
	If $ 0 < r \leq \eta $ and $ v \in \mathbf{B}(0,r) \cap \Tan_{*}(B,a) \sim \{0\} $ then we choose $ i = 1 , \ldots , l $ such that $ |(v/|v|) - v_{i}| \leq \epsilon $ and, since $ \Lip \bm{\delta}_{B} \leq 1 $, we conclude that $ \bm{\delta}_{B}(a+v) \leq 2\epsilon |v| $.
	
	\textit{Proof of \eqref{basic remark on pt tangent cone:3}}. For the necessity, suppose $ M $ is as in \ref{definition of pointwise differentiability} when $ k =1 $ and $ \alpha = 0 $, observe that $ \Tan(M,a) = \Tan^{*}(B,a) $ by \cite[3.4]{snulmenn:sets.v1} and $ \Tan(M,a) \subseteq \Tan_{*}(B,a) $ because
	\begin{equation*}
		\lim_{\Tan(M,a) \ni v \to 0} |v|^{-1}\bm{\delta}_{M}(a+v) = 0.
	\end{equation*}
	For the sufficiency let $ M = \{ a+v: v \in \Tan(B,a) \} $ and, since $ a \in \textrm{Clos}\,B $, one verifies that
	\begin{eqnarray*}
		& \sup\{|\bm{\delta}_{B}(x)-\bm{\delta}_{M}(x)| : x \in \mathbf{B}(a,r)\} \leq & \\
		&  \leq \sup (\{ \bm{\delta}_{B}(x) : x \in \mathbf{B}(a,2r) \cap M \} \cup \{ \bm{\delta}_{M}(x) : x \in \mathbf{B}(a,2r)\cap B \}), &
	\end{eqnarray*}
	Therefore the conclusion comes from \eqref{basic remark on pt tangent cone:1} and \eqref{basic remark on pt tangent cone:2}.
\end{Proof}

\begin{Remark}
	Compare \ref{basic remark on pt tangent cone}\eqref{basic remark on pt tangent cone:3} with the analogous result for approximate differentiability in \ref{alternative characterization of approximately differentiability}. Moreover \ref{basic remark on pt tangent cone}\eqref{basic remark on pt tangent cone:3} is a restatement of \cite[3.19]{snulmenn:sets.v1}.
\end{Remark}

\begin{Remark}\label{approx and point cones for submanifolds}
	If $ M $ is an $ m $ dimensional submanifold of class $ 1 $ of $ \Real{n} $ then, by \ref{basic remark on pt tangent cone}\eqref{basic remark on pt tangent cone:3}, \ref{alternative characterization of approximately differentiability} and \ref{inclusion of tangent cones}, one may readily infer that 
	\begin{equation*}
		\Tan(M,a) = \Tan^{m}(\Haus{m}\restrict M,a) \quad \textrm{for every $ a \in M $.}
	\end{equation*}
\end{Remark}

\begin{Theorem}\label{ap diff and pt diff}
	Let $ 1 \leq m \leq n $ and $ k \geq 1 $ be integers, $ 0 \leq \alpha \leq 1 $, \mbox{$ A \subseteq \Real{n} $} and \mbox{$ a \in \Real{n} $}. Suppose $ A $ is approximately differentiable of order $ (k,\alpha) $ at $ a $ and $ m = \dim \ap \Tan(A,a) $.
	
	Then there exists $ B \subseteq A $ pointwise differentiable of order $ (k,\alpha) $ at $ a $ \mbox{such that}
	\begin{eqnarray*}
		& \Hdensity{m}{A \sim B}{a} = 0, &\\
		& \ap \Tan(A,a) = \Tan(B,a) = \Tan^{m}(\Haus{m}\restrict B,a), &\\
		& \pt \Der^{i}B(a,\Tan(B,a)) = \ap \Der^{i}A(a) \quad \textrm{for $ i = 2 , \ldots , k $}. &
	\end{eqnarray*}
\end{Theorem}

\begin{Proof}
	Assume $ a = 0 $ and suppose $ T = \ap \Tan(A,0) $, $ P : T \rightarrow T^{\perp} $ is defined by
	\begin{equation*}
		\textstyle P(\chi) = \sum_{j=2}^{k} \langle \chi^{j}/j!, \ap \Der^{j}A(0) \rangle  \quad \textrm{for $ \chi \in T $}
	\end{equation*}
	and $ \Gamma = \sup\{1,\sum_{j=2}^{k}\| \ap \Der^{j}A(0) \| / j!\} $. In particular if $ k = 1 $ then $ P = 0 $ and $ \Gamma = 1 $. By \ref{distance vs vertical distance} and \ref{alternative characterization of approximate differentiability for sets} we infer that
	\begin{eqnarray*}
		& \Hdensity{m}{A \sim \mathbf{X}_{k}(0,T,P,\epsilon)}{a} =0 \quad \textrm{for every $ \epsilon > 0 $ if $ \alpha = 0 $}, &\\
		& \Hdensity{m}{A \sim \mathbf{X}_{k,\alpha}(0,T,P,\lambda)}{a} =0 \quad \textrm{for some $ 0 \leq \lambda < \infty $ if $ \alpha > 0 $.} &
	\end{eqnarray*}
	We fix $ 0 \leq \lambda < \infty $ as above if $ \alpha > 0 $. We define, for every integer $ i \geq 1 $,
	\begin{center}
		$  A_{i} = A \cap \mathbf{X}_{k}(0,T,P,(2i)^{-1}) $ if $ \alpha = 0, \quad  A_{i} = A \cap \mathbf{X}_{k,\alpha}(0,T,P,\lambda) $ if $ \alpha > 0 $.
	\end{center}
	Let $ Q_{r} = \Real{n}\cap \{ z: |T_{\natural}^{\perp}(z)| \leq r,\; |T_{\natural}(z)|\leq r \} $ for $ 0 < r < \infty $. For every integer $ i \geq 1 $ let $ \delta_{i} > 0 $ be such that 
	\begin{equation*}
		\Haus{m}( A \cap Q_{r} \sim A_{i} ) \leq 2^{-i}\,\bm{\alpha}(m)\,r^{m} \quad \textrm{whenever $ 0 < r \leq \delta_{i} $}
	\end{equation*}
	and we assume $ \delta_{i+1} < \delta_{i} $, $ \delta_{i} \to 0 $ as $ i \to \infty $,
	\begin{equation*}
		\delta_{1} \leq (2\Gamma)^{-1} \;\; \textrm{if $ \alpha = 0 $,} \quad \delta_{1} \leq (\lambda + \Gamma)^{-1/\alpha} \;\; \textrm{if $ \alpha > 0 $.}
	\end{equation*}
	We define, for every integer $ i \geq 1 $,
	\begin{equation*}
		C_{i} = T^{-1}_{\natural}[\mathbf{B}(0,\delta_{i}) \sim \mathbf{B}(0,\delta_{i+1})], \quad \textstyle B = \bigcup_{j=1}^{\infty} A_{j} \cap C_{j}.
	\end{equation*}
	Observe that $ B \subseteq \mathbf{X}(0,T,1) $ and
	\begin{equation*}
		(Q_{\delta_{j}} \sim Q_{\delta_{j+1}}) \cap \mathbf{X}(0,T,1) \sim C_{j} = \varnothing \quad \textrm{whenever $ j \geq 1 $.}
	\end{equation*}
	We can prove now that $ \Hdensity{m}{A \sim B}{0} = 0 $. In fact, by \ref{remark on approximate tangent cone I} and \ref{basic characterization} we infer $ \Hdensity{m}{A \sim \mathbf{X}(0,T,1)}{0} = 0 $. Moreover, if $ 0 < r \leq \delta_{1} $ and $ i \geq 1 $ are such that $ \delta_{i+1} < r \leq \delta_{i} $ then
	\begin{eqnarray*}
		&  Q_{r} \cap A \cap \mathbf{X}(0,T,1) \sim B \subseteq (Q_{r} \cap A \sim A_{i}) \cup \bigcup_{j=i+1}^{\infty} Q_{\delta_{j}} \cap A \sim A_{j}, &\\
		& \Haus{m}( Q_{r} \cap A \cap \mathbf{X}(0,T,1) \sim B ) \leq \bm{\alpha}(m)r^{m} \sum_{j=i}^{\infty} 2^{-j}. &
	\end{eqnarray*}
	Since this implies $ 0 \in \Clos B $ by \ref{basic characterization}, it follows that
	\begin{equation*}
		\lim_{r \to 0}r^{-k}\sup\{|P(T_{\natural}(z))-T^{\perp}_{\natural}(z)|: z \in B  \cap T^{-1}_{\natural} [\mathbf{B}(0,r)]  \} =0 \quad \textrm{if $ \alpha = 0 $}, 
	\end{equation*}
	\begin{equation*}
		\limsup_{r \to 0}r^{-k-\alpha}\sup\{|P(T_{\natural}(z))-T^{\perp}_{\natural}(z)|: z \in B  \cap T^{-1}_{\natural} [\mathbf{B}(0,r)]  \} \leq \lambda\quad \textrm{if $ \alpha > 0 $}.
	\end{equation*}
	In particular, $ \Tan^{*}(B,0) \subseteq T $. By \ref{Density and inclusion of approximate tangent cones} and \ref{basic characterization} we get that 
	\begin{equation*}
		T =\Tan^{m}_{*}(\Haus{m}\restrict A,0) \subseteq \Tan^{m}_{*}(\Haus{m}\restrict B,0).
	\end{equation*}
	Therefore $ B $ is pointwise differentiable of order $ 1 $ at $ a $ with $ T = \Tan(B,0) = \Tan^{m}(\Haus{m}\restrict B,0) $ by \ref{inclusion of tangent cones} and \ref{basic remark on pt tangent cone}\eqref{basic remark on pt tangent cone:3}. Moreover, since $ \Tan(\gr P, 0) = T $, we can use \ref{distance vs vertical distance} to check that the conditions in \ref{definition of pointwise differentiability} hold with $ M $ replaced by $ \gr P $. Therefore, by \cite[3.12]{snulmenn:sets.v1} and \ref{definition of approximate differentials}, we conclude that 
	\begin{equation*}
		\pt \Der^{i}B(0,T) = \ap \Der^{i}A(0) \quad \textrm{for $ i = 2 , \ldots , k $}.
	\end{equation*}

\end{Proof}

\begin{Theorem}\label{touching by balls and estimate ssf}
	Let $ A \subseteq \Real{n} $, $ a \in \Real{n} $, $ \nu \in \mathbf{S}^{n-1} $, $ 0 < r < \infty $ and suppose 
	\begin{equation*}
		\mathbf{U}(a+r\nu,r) \cap A = \varnothing.
	\end{equation*}
	
	Then the following three statements hold.
	\begin{enumerate}
		\item \label{touching by balls and estimate ssf:smooth case} If $ A $ is a submanifold of class $ 2 $ and $ a \in A $ then 
		\begin{equation*}
			\mathbf{b}_{A}(a)(v,v) \bullet \nu \leq r^{-1}|v|^{2} \quad \textrm{whenever $ v \in \Tan(A,a) $.}
		\end{equation*}
		\item \label{touching by balls and estimate ssf:pt case} If $ A $ is pointwise differentiable of order $ 2 $ at $ a $ then 
		\begin{equation*}
			\pt\Der^{2}A(a, \Tan(A,a))(v,v) \bullet \nu \leq r^{-1}|v|^{2} \quad \textrm{whenever $ v \in \Tan(A,a) $.}
		\end{equation*}
		\item \label{touching by balls and estimate ssf:ap case} If $ A $ is approximately differentiable of order $ 2 $ at $ a $ then 
		\begin{equation*}
			\ap\Der^{2}A(a)(v,v) \bullet \nu \leq r^{-1}|v|^{2} \quad \textrm{whenever $ v \in \ap \Tan(A,a) $}.
		\end{equation*}
	\end{enumerate}
\end{Theorem}

\begin{Proof}
	Assume $ a = 0 $. Observe that $ \nu \in \Nor^{*}(A,0) $.

	The statement in \eqref{touching by balls and estimate ssf:smooth case} is classical. We give a proof here for completeness. \mbox{If $ T = \Tan(A,0) $} then there exist a function $ f : T \rightarrow T^{\perp} $ of class $ 2 $ and an open neighbourhood $ U $ of $ 0 \in \Real{n} $ such that $ \Der f (0) = 0 $, $ T_{\natural}[U] = T_{\natural}[U \cap A] $ and $  A \cap U  = \{ \chi + f(\chi) : \chi \in T_{\natural}[U] \} $. Since for every $ \chi \in T_{\natural}[U] $  
	\begin{equation*}
		| \chi + f(\chi) - r\nu| \geq r, \quad 2r\,f(\chi) \bullet \nu \leq |\chi|^{2} + |f(\chi)|^{2},
	\end{equation*}
	we conclude that $ \Der^{2}f(0)(v,v) \bullet \nu \leq r^{-1}|v|^{2} $ for every $ v \in \Tan(A,0) $ and, since $ \mathbf{b}_{A}(0) = \Der^{2}f(0) $, the statement in \eqref{touching by balls and estimate ssf:smooth case} follows.
	
	The statement in \eqref{touching by balls and estimate ssf:pt case} is mainly a consequence of \cite[3.18]{snulmenn:sets.v1}. In fact suppose $ T = \Tan(A,0) $, $ P : T \rightarrow T^{\perp} $ is the homogeneous polynomial function of degree $ 2 $ such that $ \pt\Der^{2}A(0,T) = \Der^{2}(P \circ T_{\natural})(0) $ (whose existence can be asserted, from instance, by \cite[3.22]{snulmenn:sets.v1}) and \mbox{$ B = \{ \chi + P(\chi) : \chi \in T \} $.} If we prove that $ \mathbf{U}(r\nu,r) \cap B = \varnothing $ then \eqref{touching by balls and estimate ssf:pt case} is a consequence of \eqref{touching by balls and estimate ssf:smooth case}. By contradiction let $ x \in B \cap \mathbf{U}(r\nu,r) $ and, by \cite[3.18]{snulmenn:sets.v1}, for every positive integer $ i $ we can select $ x_{i} \in A $ such that 
	\begin{equation*}
		|iT_{\natural}(x_{i}) + i^{2}T_{\natural}^{\perp}(x_{i}) - x | \to 0 \quad \textrm{as $ i \to \infty $}.
	\end{equation*}
	Since $ | x_{i} - r\nu| \geq r $ for every $ i \geq 1 $, we get
	\begin{flalign*}
		&	|iT_{\natural}(x_{i}) + i^{2}T^{\perp}_{\natural}(x_{i}) -r\nu|^{2} \\
		&  \quad = i^{2} |x_{i}-r\nu|^{2} + (i^{4}-i^{2})|T_{\natural}^{\perp}(x_{i})|^{2} + r^{2}- i^{2}r^{2} \\
		& \quad \geq (i^{4}-i^{2})|T_{\natural}^{\perp}(x_{i})|^{2} + r^{2} \quad \textrm{for $ i \geq 1 $;}
	\end{flalign*}
	yet $ | iT_{\natural}(x_{i}) + i^{2}T_{\natural}^{\perp}(x_{i}) - r\nu| < r $ for $ i $ large. This is a contradiction.
	
	Finally \eqref{touching by balls and estimate ssf:ap case} is a consequence of \eqref{touching by balls and estimate ssf:pt case} and \ref{ap diff and pt diff}.
\end{Proof}

\section{Rectifiability and Borel measurability}\label{section: Rectifiability and Borel measurability}

\begin{Lemma}\label{Borel measurability 1}
	Let $ 1 \leq m \leq n $  and $ k \geq 1 $ be integers, $ 0 \leq \alpha \leq 1 $, $ \gamma = k + \alpha $ and $ A \subseteq \Real{n} $. Let $ Y $ be the set of
	\begin{equation*}
		(a,T,\phi_{0}, \ldots, \phi_{k}) \in \Real{n} \times \mathbf{G}(n,m) \times \prod_{i=0}^{k}\textstyle \bigodot^{i}(\Real{n}, \Real{n})
	\end{equation*} 
	such that $ \phi_{0} = T_{\natural}^{\perp}(a) $ and 
	\begin{equation*}
		\lim_{r \to 0}\frac{\Haus{m}\left( A \cap \mathbf{U}(a,r)\cap \{ z: | T^{\perp}_{\natural}(z) - \sum_{j=0}^{k}\langle T_{\natural}(z-a)^{j}/j!,\phi_{j} \rangle | > \lambda \,r^{\gamma}  \} \right) }{\bm{\alpha}(m)\,r^{m}} =0
	\end{equation*}
	for every $ \lambda > 0 $ [for some $ 0 \leq \lambda < \infty $].
	
	Then $ Y $ is a Borel subset of $ \Real{n} \times \mathbf{G}(n,m) \times \prod_{i=0}^{k}\textstyle \bigodot^{i}(\Real{n}, \Real{n}) $.
\end{Lemma}

\begin{Proof}
	Let $ Z = \Real{n} \times \mathbf{G}(n,m) \times \prod_{j=0}^{k}\textstyle \bigodot^{j}(\Real{n}, \Real{n}) $. If $ 0 < \lambda < \infty $, $ i \geq 1 $ is an integer and $ 0 < r < \infty $, we define $ W_{\lambda,i,r} $ to be the set of $ (a,T,\phi_{0}, \ldots ,\phi_{k}) \in Z $ such that $ \phi_{0} = T^{\perp}_{\natural}(a) $ and 
	\begin{equation*}
		\textstyle \Haus{m}\left(A \cap \mathbf{U}(a,r) \cap \left\{ z: | T^{\perp}_{\natural}(z) - \sum_{l=0}^{k}\langle T_{\natural}(z-a)^{l}/l!,\phi_{l} \rangle | > \lambda \,r^{\gamma}  \right\} \right) \leq i^{-1}\,r^{m}.
	\end{equation*}
	Then $ W_{\lambda,i,r} $ is a closed subset of $ Z $. In fact if $ (a_{j},T_{j},\phi_{0,j}, \ldots ,\phi_{k,j}) \in W_{\lambda,i,r} $, $ j \geq 1 $, is a sequence converging to $ (a,T,\phi_{0}, \ldots ,\phi_{k}) \in Z $ as $ j \to \infty $, we define 
	\begin{eqnarray*}
		& P_{j}(\chi) = \sum_{l=0}^{k}\langle (\chi-T_{j\,\natural}(a_{j}))^{l}/l!,\phi_{l,j} \rangle \quad \textrm{for $ \chi \in T_{j} $ and $ j \geq 1 $}, &\\
		& P(\chi) = \sum_{l=0}^{k}\langle (\chi-T_{\natural}(a))^{l}/l!,\phi_{l} \rangle \quad \textrm{for $ \chi \in T $}, &
	\end{eqnarray*}
	and we observe that $ P_{j}(T_{j\,\natural}(z)) \to P(T_{\natural}(z)) $ as $ j \to \infty $, whenever $ z \in \Real{n} $. Let
	\begin{eqnarray*}
		& B_{j} = A \cap \mathbf{U}(a_{j},r)\cap \{ z: | T^{\perp}_{j\,\natural}(z) - P_{j}(T_{j\,\natural}(z)) | > \lambda \,r^{\gamma}  \}, &\\
		& B = A \cap \mathbf{U}(a,r)\cap \{ z: | T^{\perp}_{\natural}(z) - P(T_{\natural}(z)) | > \lambda \,r^{\gamma}  \} &
	\end{eqnarray*}
	and observe that
	\begin{equation*}
		B \subseteq \bigcup_{j=1}^{\infty}\bigcap_{h=j}^{\infty}B_{h}.
	\end{equation*}
	Therefore, by \cite[2.1.5(1)]{MR0257325}, we conclude that
	\begin{equation*}
		\Haus{m}(B) \leq \lim_{j \to \infty} \Haus{m}( \bigcap_{h=j}^{\infty}B_{h})\leq
		\liminf_{j \to \infty} \Haus{m}(B_{j}) \leq i^{-1} r^{m},
	\end{equation*}
	$ (a,T,\phi_{0}, \ldots, \phi_{k}) \in W_{\lambda,i,r} $ and $ W_{\lambda,i,r} $ is closed. 
	
	Henceforth $ Y $ is a Borel set because
	\begin{equation*}
		\textstyle Y = \bigcap_{l=1}^{\infty}\bigcap_{i=1}^{\infty}\bigcup_{j=1}^{\infty}\bigcap\{ W_{l^{-1},i,r}: 0 < r \leq j^{-1}\}, 
	\end{equation*}
	\begin{equation*}
		\textstyle \big[Y = \bigcup_{l=1}^{\infty}\bigcap_{i=1}^{\infty}\bigcup_{j=1}^{\infty}\bigcap\{ W_{l,i,r}: 0 < r \leq j^{-1}\}\big].
	\end{equation*}
\end{Proof}

\begin{Lemma}\label{Borel measurability 2}
	Suppose $ 1 \leq m \leq n $ are integers, $ A \subseteq \Real{n} $ and $ \tau_{a}(x) = x-a $ whenever $ a,x \in \Real{n} $. Let $ Y $ be the set of $ (a,T) \in \Real{n} \times \mathbf{G}(n,m) $ such that for every $ \epsilon > 0 $ there exist $ \eta > 0 $ and $ \rho > 0 $ such that 
	\begin{equation*}
		\Haus{m}(\mathbf{C}(T,z,\epsilon r, \epsilon r) \cap \tau_{a}[A] ) \geq \eta \bm{\alpha}(m) r^{m}
	\end{equation*}
	for every $ 0 < r \leq \rho $ and for every $ z \in T \cap \mathbf{B}(0,r) $.
	
	Then $ Y $ is a Borel subset of $ \Real{n} \times \mathbf{G}(n,m) $.
	
\end{Lemma}

\begin{Proof}
	We prove that $ (\Real{n} \times \mathbf{G}(n,m))  \sim Y $ is a Borel subset of $ \Real{n} \times \mathbf{G}(n,m) $. For every $ \epsilon> 0 $, $ \eta > 0 $ and $ 0 < \rho_{2} < \rho_{1} $ suppose $ W_{\epsilon,\eta,\rho_{1},\rho_{2}} $ is the set of $ (a,T) \in \Real{n} \times \mathbf{G}(n,m) $ such that
	\begin{equation*}
		\Haus{m}(\mathbf{C}(T,z,\epsilon r , \epsilon r ) \cap \tau_{a}[A]) \leq \eta \bm{\alpha}(m) r^{m}
	\end{equation*}
	for some $ z \in \mathbf{B}(0,r) \cap T $ and some $ \rho_{2} \leq r \leq \rho_{1} $. We prove that $ W_{\epsilon,\eta, \rho_{1}, \rho_{2}} $ is a closed subset of $ \Real{n}\times \mathbf{G}(n,m)  $. Suppose $ (a_{j},T_{j}) \in W_{\epsilon, \eta,\rho_{1}, \rho_{2}} $, $ j \geq 1 $, is a sequence converging to $ (a,T) \in \Real{n}\times \mathbf{G}(n,m) $ as $ j \to \infty $. Therefore there exist sequences $ \rho_{2} \leq r_{j} \leq \rho_{1} $ and $ z_{j} \in T_{j} \cap \mathbf{B}(0,r_{j}) $, for $ j \geq 1 $, such that 
	\begin{equation*}
		\Haus{m}(\mathbf{C}(T_{j},z_{j},\epsilon r_{j}, \epsilon r_{j}) \cap \tau_{a_{j}}[A]) \leq \eta \bm{\alpha}(m) r_{j}^{m} \quad \textrm{for every $ j \geq 1 $.}
	\end{equation*}
	Then there exist $ z \in \Real{n} $ and $ r \in \Real{} $ such that, possibly passing to a subsequence, $ z_{j} \to z $ and $ r_{j} \to r $ as $ j \to \infty $. Observe that $ z \in \mathbf{B}(0,r) \cap T $ and $ \rho_{2} \leq r \leq \rho_{1} $. For each $ j \geq 1 $ we define
	\begin{equation*}
		B_{j}= \mathbf{C}(T_{j},z_{j},\epsilon r_{j}, \epsilon r_{j}) \cap \tau_{a_{j}}[A], \quad
		B =\mathbf{C}(T,z,\epsilon r, \epsilon r) \cap \tau_{a}[A],
	\end{equation*}
and one may easily verify that 
\begin{equation*}
 B \subseteq \bigcup_{h=1}\bigcap_{k=h}^{\infty}\tau_{a-a_{k}}[B_{k}].
\end{equation*}
Now we can use \cite[2.1.5(1)]{MR0257325} to conclude that 
	\begin{equation*}
		\Haus{m}(B) \leq \liminf_{h \to \infty} \Haus{m}\big(\tau_{a-a_{h}}[B_{h}]\big) \leq \bm{\alpha}(m) \eta r^{m}.
	\end{equation*}
	Therefore $ (a,T) \in W_{\epsilon,\eta,\rho_{1}, \rho_{2}} $ and $ W_{\epsilon,\eta,\rho_{1}, \rho_{2}} $ is a closed subset of $ \Real{n} \times \mathbf{G}(n,m) $. If $ E \subseteq \Real{} $ is a countable set such that $ \inf E = 0 \notin E $ then it is not difficult to see that 
	\begin{equation*}
		\big(\Real{n} \times \mathbf{G}(n,m)\big) \sim Y = \bigcup_{\epsilon \in E} \bigcap_{\eta \in E} \bigcap_{\rho_{1} \in E} \bigcup_{\rho_{2} \in E} W_{\epsilon, \eta, \rho_{1}, \rho_{2}}.
	\end{equation*}
\end{Proof}

\begin{Definition}
	A measure $ \phi $ over $ X $ is called $\sigma$ finite if there exists a sequence $ X_{i} $ such that $ \phi(X_{i}) < \infty $ for every $ i \geq 1 $ and $ X = \bigcup_{i=1}^{\infty}X_{i} $.
\end{Definition}

\begin{Theorem}\label{criterion for higher order rectifiability}
	Suppose $ 1 \leq m \leq n $ and $ k \geq 1 $ are integers, $ 0 \leq \alpha \leq 1 $, $ A \subseteq \Real{n} $ such that $ \Haus{m} \restrict A $ is $\sigma$ finite and for every $ a \in A $ there exists an $ m $ dimensional submanifold $ B \subseteq \Real{n} $ of class $ (k,\alpha) $ such that $ a \in B $ and the following condition $ (\ast) $ is satisfied. \mbox{For every $ \epsilon > 0 $}
	\begin{equation*}
		\lim_{r \to 0} \frac{ \Haus{m}\left( A \cap \mathbf{B}(a,r) \cap \{ z: \bm{\delta}_{B}(z) > \epsilon\,r^{k} \} \right) }{\bm{\alpha}(m)r^{m}  }  =0 
	\end{equation*}
	and, if $ \alpha > 0 $, there exists $ 0 \leq \lambda < \infty $ such that
	\begin{equation*}
		\lim_{r \to 0} \frac{ \Haus{m}\left( A \cap \mathbf{B}(a,r) \cap \{ z: \bm{\delta}_{B}(z) > \lambda \,r^{k+\alpha} \} \right) }{\bm{\alpha}(m)r^{m}  }  = 0.
	\end{equation*}
	
	Then $ A $ is countably $ \rect{m} $ rectifiable of class $ (k,\alpha) $.
\end{Theorem}

\begin{Proof}

Clearly, we can assume $ \Haus{m}(A)< \infty $. If $ X $ is the set of points $ a \in \Real{n} $ such that there exists an $ m $ dimensional submanifold $ B $ of class $ (k,\alpha) $ such that $ a \in B $ and $ (\ast) $ is satisfied, then $ X $ is an $ \Haus{m} $ measurable subset of $ \Real{n} $ by \ref{distance vs vertical distance}, \ref{Borel measurability 1} and \cite[2.2.13]{MR0257325}. If $ E \subseteq \Real{n} $ is an $ \Haus{m} $ hull of $ A $, by \cite[2.1.5(2)]{MR0257325}, then $ A \subseteq E \cap X $ and $ E \cap X $ is a $ \Haus{m} $ measurable subset of $ \Real{n} $ that satifies the same hypothesis $ A $ does. Therefore we can assume $ A $ to be $ \Haus{m} $ measurable.
	
	If $ a \in A $ and $ B $ is an $ m $ dimensional submanifold of class $ 1 $ such that $ a \in B $ and $ (\ast) $ is satisfied then, by \ref{distance vs vertical distance} and \ref{remark on approximate tangent cone I}, we get that
	\begin{equation*}
	\Tan^{*m}(\Haus{m}\restrict A,a) \subseteq \Tan(B,a), \quad	\Hdensity{m}{A \sim \mathbf{X}(a,\Tan(B,a),\epsilon)}{a} = 0, 
	\end{equation*}
	for every $ \epsilon > 0 $. Therefore by \cite[2.10.19(2), 3.3.17, 3.2.29]{MR0257325} we conclude that $ A $ is $ \rect{m} $ rectifiable of class $ 1 $.
	
	Let $ S \in \mathbf{G}(n,m) $, let $ U \subseteq S $ be relatively open, let $ f : U \rightarrow S^{\perp} $ be a function of \mbox{class $ 1 $}, $ M = \{  \chi + f(\chi) : \chi \in U  \} $, $ \Lip f < \infty $ and $ \Haus{m}(M)< \infty $. We prove that $ A \cap M $ is $ \rect{m} $ rectifiable of class $ (k,\alpha) $. This evidently implies that $ A $ is $ \rect{m} $ rectifiable of class $ (k,\alpha) $. Let $ Y $ be the set of points $ a \in A \cap M $ such that $ \Hdensity{m}{M \sim A}{a}=0 $. We use \ref{inclusion of tangent cones}, \ref{approx and point cones for submanifolds} and \ref{Density and inclusion of approximate tangent cones} to conclude that
	\begin{equation*}
		\Tan^{*m}(\Haus{m}\restrict A \cap M,a) = \Tan(M,a) \quad \textrm{for every $ a \in Y $}.
	\end{equation*}
	By \cite[2.10.19(4)]{MR0257325} we have $ \Haus{m}( A \cap M \sim Y)=0 $. Let 
	\begin{equation*}
		C = S \cap \{  \chi : \chi + f(\chi) \in A \cap M  \}, \quad D = S \cap \{ \chi: \chi + f(\chi)\in Y   \},
	\end{equation*}
	we observe that $ \Haus{m}( C \sim D)=0 $ and
	\begin{equation*}
		\Hdensity{m}{S\sim C}{\chi} =0 \quad \textrm{for every $ \chi \in D $.}
	\end{equation*}
	
	Let $ \chi \in D $, $ a = \chi + f(\chi) $ and suppose $ B $ is an $ m $ dimensional submanifold of class $ (k,\alpha) $ such that $ a \in B $ and $ (\ast) $ is satisfied. Since $ \Tan(B,a) \cap S^{\perp} = \{0\} $, there exist a function $ g : S \rightarrow S^{\perp} $ of class $ (k,\alpha) $ and an open neighbourhood $ V $ of $ a $ such that $ B \cap V = \{ \zeta + g(\zeta) : \zeta \in S  \} \cap V $. Therefore, by \ref{distance vs vertical distance},
	\begin{equation*}
		\lim_{r \to 0}\frac{ \Haus{m}\left( A \cap \mathbf{B}(a,r) \cap \{ z: | g(S_{\natural}(z)) - S^{\perp}_{\natural}(z)| > \epsilon \,r^{k} \} \right) }{\bm{\alpha}(m)r^{m}}  = 0 
	\end{equation*}
	for every $ \epsilon > 0 $ and, if $ \alpha > 0 $, there exists $ 0 \leq \lambda < \infty $ such that
	\begin{equation*}
		\lim_{r \to 0}\frac{ \Haus{m}\left( A \cap \mathbf{B}(a,r) \cap \{ z: | g(S_{\natural}(z)) - S^{\perp}_{\natural}(z)| > \lambda \,r^{k+\alpha} \} \right) }{\bm{\alpha}(m)r^{m}}  = 0.
	\end{equation*}
	Let $ P : S \rightarrow S^{\perp} $ be the $ k $ jet of $ g $ at $ \chi $. If $ \epsilon > 0 $ then, possibly replacing $ \lambda $ by a larger number if $ \alpha > 0 $, we can choose $ \rho > 0 $ such that 
	\begin{equation*}
		|g(\zeta) - P(\zeta)| \leq \lambda\,r^{k+\alpha} \;\; \textrm{if $ \alpha > 0 $}, \quad |g(\zeta) - P(\zeta)| \leq \epsilon\,r^{k} \;\; \textrm{if $ \alpha = 0 $,}
	\end{equation*}
	for every $ \zeta \in \mathbf{B}(\chi,r) $ and $ 0 < r \leq \rho $. Let $ \Gamma =(1+ (\Lip f)^{2})^{1/2} $,  $ \gamma = \lambda + \Gamma^{k+\alpha}\,\lambda $ if $ \alpha > 0 $ and observe that, whenever $ 0 < r \leq \rho $,
	\begin{flalign*}
		&	C \cap \mathbf{B}(\chi,r) \cap \{ \zeta: |f(\zeta) - P(\zeta)| > \gamma\,r^{k+\alpha} \} \\
		& \quad  \subseteq S_{\natural} [ A \cap M \cap \mathbf{B}(a,\Gamma\,r) \cap \{ z: |S^{\perp}_{\natural}(z) - g(S_{\natural}(z))| >\lambda\,\Gamma^{k+\alpha}\,r^{k+\alpha}   \}] \quad \textrm{if $ \alpha > 0 $}, \\
		& C \cap \mathbf{B}(\chi,r) \cap \{ \zeta: |f(\zeta) - P(\zeta)| > 2\epsilon\,r^{k} \} \\
		& \quad  \subseteq S_{\natural} [ A \cap M \cap \mathbf{B}(a,\Gamma\,r) \cap \{ z: |S^{\perp}_{\natural}(z) - g(S_{\natural}(z))| >\epsilon\,\,r^{k}   \}] \quad \textrm{if $ \alpha = 0 $}.
	\end{flalign*}
	Since $ \chi $ is arbitrarily chosen in $ D $, we infer by \ref{alternative characterization of approximate differentiability for functions} and \ref{rectifiability and approximate differentiability for functions} that there exist countably many functions $ g_{j} : S \rightarrow S^{\perp} $ of class $ (k,\alpha) $ such that
	\begin{equation*}
		\Haus{m}\left(C \sim \textstyle \bigcup_{j=1}^{\infty}\{ \zeta : f(\zeta) = g_{j}(\zeta)  \} \right) =0,
	\end{equation*}
	whence $ A \cap M $ is $ \rect{m} $ rectifiable of class $ (k,\alpha) $. 
\end{Proof}

\begin{Theorem}\label{approx differentials are Borel maps}
	Suppose $ n \geq 1 $, $ k \geq 1 $ are integers, $ 0 \leq \alpha \leq 1 $, $ A \subseteq \Real{n} $ and $ X $ is the set of $ a \in \Real{n} $ where $ A $ is approximately differentiable of order $(k,\alpha)$.
	
	Then $ X $ is a Borel subset of $ \Real{n} $. Moreover $ \ap \Tan(A,\cdot)_{\natural} $ is a Borel map whose domain is a Borel subset of $ \Real{n} $ and the same conclusion is true for $ \ap \Der^{j} A $ for every $ j \geq 2 $. 
\end{Theorem}

\begin{Proof}
	First apply \ref{Borel measurability 1} and \ref{Borel measurability 2} to see that the set $ Z $ of 
\begin{equation*}
 (a,T, \phi_{0}, \ldots , \phi_{k}) \in \bigcup_{m=1}^{n}\big[\Real{n} \times \mathbf{G}(n,m) \times \textstyle \prod_{i=0}^{k}\textstyle \bigodot^{i}(\Real{n}, \Real{n})\big]
\end{equation*}
such that $ \phi_{1}= 0 $ and the conditions listed in \ref{Borel measurability 1} and \ref{Borel measurability 2} are satisfied for every $ \lambda > 0 $ if $ \alpha = 0 $ and for some $ 0 \leq \lambda < \infty $ if $ \alpha > 0 $, is a Borel set; then use \ref{distance vs vertical distance} and \ref{basic characterization} to conclude that $ Z $ is the graph of a function $ f $ mapping $ X $ into $ \bigcup_{m=1}^{n}[\mathbf{G}(n,m) \times \prod_{i=0}^{k}\textstyle \bigodot^{i}(\Real{n}, \Real{n})] $; finally apply \cite[4.1]{snulmenn:sets.v1} to infer that $ X $ is a Borel subset of $ \Real{n} $ and $ f $ is a Borel function. In case $ \alpha = 0 $, this clearly proves the second part of the conclusion.
\end{Proof}

\begin{Theorem}\label{the set of approximate differentiable points is rectifiable}
	Suppose $ 1 \leq m \leq n $ and $ k \geq 1 $ are integers, $ 0 \leq \alpha \leq 1 $, $ A \subseteq \Real{n} $ is $ \Haus{m} $ measurable with $ \Haus{m}(A)< \infty $ and $ X $ is the set of $ a \in \Real{n} $ such that $ A $ is approximately differentiable of order $ (k,\alpha) $ at $ a $ with $ \dim \ap \Tan(A,a) = m $.
	
	Then $ X $ is $ \rect{m} $ rectifiable of class $ (k,\alpha) $ and $ \Haus{m}(X \sim A) = 0 $.
\end{Theorem}

\begin{Proof}
	Apply \ref{criterion for higher order rectifiability} to get that $ A \cap X $ is $ \rect{m} $ rectifiable of class $ (k,\alpha) $. Since $ \infHdensity{m}{A}{x} > 0 $ for every $ x \in X $ by \ref{basic characterization}, we infer that $\Haus{m}(X \sim A) = 0$ by \cite[2.10.19(4)]{MR0257325}.
\end{Proof}

\begin{Remark}
	The pattern of this section follows \cite[\S 5]{snulmenn:sets.v1}.  
\end{Remark}

\medskip

{\small \noindent Max Planck Institute for Gravitational Physics (Albert
Einstein Institute) \newline Am M{\"u}hlen\-berg 1, \newline D-14476 Golm,
Germany \newline \texttt{mario.santilli@aei.mpg.de}}

\vspace{4mm}

{\small \noindent University of Potsdam, Institute for Mathematics, \newline
OT Golm, Karl-Liebknecht-Stra{\ss}e 24-25, \newline
	D-14476 Potsdam, Germany}

\begin{thebibliography}{{Men}16}
	
	\bibitem[Alb94]{MR1384392}
	Giovanni Alberti.
	\newblock On the structure of singular sets of convex functions.
	\newblock {\em Calc. Var. Partial Differential Equations}, 2(1):17--27, 1994.
	
	\bibitem[All72]{MR0307015}
	William~K. Allard.
	\newblock On the first variation of a varifold.
	\newblock {\em Ann. of Math. (2)}, 95:417--491, 1972.
	
	\bibitem[AS94]{MR1285779}
	Gabriele Anzellotti and Raul Serapioni.
	\newblock {$\mathscr C^k$}-rectifiable sets.
	\newblock {\em J. Reine Angew. Math.}, 453:1--20, 1994.
	
	\bibitem[BHS05]{MR2183045}
	Bogdan Bojarski, Piotr Haj{\l}asz, and Pawe{\l} Strzelecki.
	\newblock Sard's theorem for mappings in {H}\"older and {S}obolev spaces.
	\newblock {\em Manuscripta Math.}, 118(3):383--397, 2005.
	
	\bibitem[Cam64]{MR0167862}
	S.~Campanato.
	\newblock Propriet\`a di una famiglia di spazi funzionali.
	\newblock {\em Ann. Scuola Norm. Sup. Pisa (3)}, 18:137--160, 1964.
	
	\bibitem[Del12]{MR2904134}
	Silvano Delladio.
	\newblock Functions of class {$C^1$} subject to a {L}egendre condition in an
	enhanced density set.
	\newblock {\em Rev. Mat. Iberoam.}, 28(1):127--140, 2012.
	
	\bibitem[Fed59]{MR0110078}
	Herbert Federer.
	\newblock Curvature measures.
	\newblock {\em Trans. Amer. Math. Soc.}, 93:418--491, 1959.
	
	\bibitem[Fed69]{MR0257325}
	Herbert Federer.
	\newblock {\em Geometric measure theory}.
	\newblock Die Grundlehren der mathematischen Wissenschaften, Band 153.
	Springer-Verlag New York Inc., New York, 1969.
	
	\bibitem[FM99]{MR1686704}
	Ilaria Fragal\`a and Carlo Mantegazza.
	\newblock On some notions of tangent space to a measure.
	\newblock {\em Proc. Roy. Soc. Edinburgh Sect. A}, 129(2):331--342, 1999.
	
	\bibitem[Isa87]{MR897693}
	N.~M. Isakov.
	\newblock On a global property of approximately differentiable functions.
	\newblock {\em Mat. Zametki}, 41(4):500--508, 620, 1987.
	
	\bibitem[Koh77]{MR0427559}
	Robert~V. Kohn.
	\newblock An example concerning approximate differentiation.
	\newblock {\em Indiana Univ. Math. J.}, 26(2):393--397, 1977.
	
	\bibitem[Kol16]{2015arXiv150600507K}
	S.~Kolasi{\'n}ski.
	\newblock {Higher order rectifiability of measures via averaged discrete
		curvatures}.
	\newblock {\em ArXiv e-prints}, pages 1--20, April 2016.
	
	\bibitem[Mat95]{MR1333890}
	Pertti Mattila.
	\newblock {\em Geometry of sets and measures in {E}uclidean spaces}, volume~44
	of {\em Cambridge Studies in Advanced Mathematics}.
	\newblock Cambridge University Press, Cambridge, 1995.
	\newblock Fractals and rectifiability.
	
	\bibitem[Men13]{zbMATH06157228}
	U.~Menne.
	\newblock {Second order rectifiability of integral varifolds of locally bounded
		first variation.}
	\newblock {\em {J. Geom. Anal.}}, 23(2):709--763, 2013.
	
	\bibitem[{Men}16]{snulmenn:sets.v1}
	U.~{Menne}.
	\newblock {Pointwise differentiability of higher order for sets}.
	\newblock {\em ArXiv e-prints}, pages 1--33, March 2016.
	
	\bibitem[MS17]{2017arXiv170309561M}
	U.~{Menne} and M.~{Santilli}.
	\newblock {A geometric second-order-rectifiable stratification for closed
		subsets of Euclidean space}.
	\newblock {\em ArXiv e-prints}, March 2017.
	
	\bibitem[Sim83]{MR756417}
	Leon Simon.
	\newblock {\em Lectures on geometric measure theory}, volume~3 of {\em
		Proceedings of the Centre for Mathematical Analysis, Australian National
		University}.
	\newblock Australian National University, Centre for Mathematical Analysis,
	Canberra, 1983.
	
	\bibitem[Whi51]{MR0043878}
	Hassler Whitney.
	\newblock On totally differentiable and smooth functions.
	\newblock {\em Pacific J. Math.}, 1:143--159, 1951.
	
\end{thebibliography}
\end{document}